\documentclass[11pt,a4paper]{article}
\usepackage{beamerarticle}

\mode<article>{\usepackage{fullpage}}
\mode<presentation>{\usetheme{Berlin}}

\usepackage{amsfonts, amsmath}
\usepackage{graphicx}

\usepackage{color}

\usepackage[english]{babel}
\usepackage{pgf}
\usepackage[utf8]{inputenc}
\usepackage{fontenc}

\newcommand{\DMCF}{$\psi$MCF }
\newcommand{\wkp}{{k_\psi^2}}

\newcommand{\com}[1]{\opt{draft}{\textcolor{red}{
$\LHD$ #1 $\RHD$\marginpar{\textcolor{red}{$\begin{lema}acksquare$}}}}}

\newcommand{\comb}[1]{\opt{draft}{\textcolor{blue}{
$\LHD$ #1 $\RHD$\marginpar{\textcolor{blue}{$\begin{lema}acksquare$}}}}}

\def\qed{\hfill {\large ${\sqcup\!\!\!\!\sqcap}$}}

\newenvironment{demo}{{\bf Proof }}
{\qed \\}

\newenvironment{prue}{{\it Proof }}
{\qed \\}

\newcommand{\grad}{{\rm grad}}
\newcommand{\dv}{{\rm div}}

\newcommand{\re}{\mathbb R}
\newcommand{\ene}{\mathbb N}

\newcommand{\<}{\left<}
\renewcommand{\(}{\left(}
\newcommand{\lb}{\label}
\newcommand{\nn}{\nonumber}

\newcommand{\fracc}{\displaystyle\frac}
\newcommand{\ds}{\displaystyle}
\renewcommand{\>}{\right>}
\renewcommand{\)}{\right)}
\newcommand{\flecha}{\longrightarrow}

\newcommand{\eps}{\ensuremath{\varepsilon}}

\def\p{\varphi}

\def\parcial#1#2{\frac{\partial #1}{\partial#2}}

\def\deri#1#2{\frac{d #1}{d#2}}

\def\sgn{{\rm sgn}}

\def\tr{{\rm tr}}
\def\vec{\overrightarrow}

\def\parcial#1#2{\fracc{\partial #1}{\partial#2}}

\def\deri#1#2{\fracc{d #1}{d#2}}

\def\ft{\frak t}

\def\tr{{\rm tr}}

\newcommand{\bde}{\begin{defi}}
\newcommand{\ede}{\end{defi}}

\newcommand{\be}{\begin{enumerate}}
\newcommand{\ee}{\end{enumerate}}

\newcommand{\ba}{\begin{array}}
\newcommand{\ea}{\end{array}}

\def\Hp{{H_\psi}}

\def\t{\theta}

\def\og{{\overline g}}

\def\oRic{{\overline Ric}}
\def\oM{{\overline M}}
\def\ona{{\overline \nabla}}

\def\oDelta{\overline{\Delta}}

\def\oO{{\overline \Omega}}

\def\wga{{\widetilde {\gamma} }}

\def\hOm{{\widehat \Omega }}

\def\htt{{\widehat {t}}}
\def\hF{{\widehat {F}}}
\def\hA{{\widehat {A}}}
\def\hC{{\widehat {C}}}

\newtheorem{defi}{\hspace{12pt} Definition}
\newtheorem{teor}{\hspace{12pt} Theorem}
\newtheorem{prop}[teor]{\hspace{12pt} Proposition}

\newtheorem{lema}[teor]{\hspace{12pt} Lemma}
\newtheorem*{lema*}{\hspace{12pt} Lemma}

\newcommand{\ben}{\begin{enumerate}}
\newcommand{\een}{\end{enumerate}}
\newcommand{\bi}{\begin{itemize}}
\newcommand{\ei}{\end{itemize}}
\newcommand{\bec}{\begin{equation}}
\newcommand{\eec}{\end{equation}}
\newcommand{\beca}{\begin{equation*}}
\newcommand{\eeca}{\end{equation*}}
\newcommand{\bal}{\begin{align}}
\newcommand{\aal}{\end{align}}
\newcommand{\bala}{\begin{align*}}
\newcommand{\aala}{\end{align*}}

\begin{document}

\title{The curve shortening problem associated to a density}

\author{ Vicente Miquel 
and Francisco Viñado-Lereu 
\thanks{Research partially supported by  the
DGI (Spain) and FEDER  project MTM2013-46961-P. and the Generalitat Valenciana Project  PROMETEOII/2014/064. The second author has been supported by a Grant of the Programa Nacional de Formaci\'on de Personal Investigador 2011 Subprograma FPI-MICINN ref: BES-2011-045388.}
}

\date{}

\maketitle

\vspace{-1cm}
\begin{abstract} In $\re^n$ with a density $e^\psi$, we study the mean curvature flow associated to the density ($\psi$-mean curvature flow or \DMCF) of a hypersurface. The main results concern with the description of the evolution under \DMCF of a closed embedded curve in the plane with a radial density, and with a statement of subconvergence to a $\psi$-minimal closed curve in a surface under some general circumstances.
\end{abstract}

\section{Introduction }\lb{In}

The mean curvature flow (MCF for short) of an immersion $F:M\flecha \oM$ of a hypersurface $M$ in a $n+1$ dimensional  Riemannian manifold $(\oM,\og)$ is the $1$-parametric family of immersions $F:M\times [0,T[\flecha \oM$  solution of the equation 

\begin{equation}
\parcial{F}{t} = \vec{H} = H N,\nn
\end{equation}
where $H$ is the mean curvature of the immersion, $N$ is a unit vector field orthogonal to the immersed hypersurface, and we have used the following {\it convention signs for the mean curvature $H$, the Weingarten map $A$ and the second fundamental form $h$}:

$A X = - \ona_XN$, $h(X,Y) = \<\ona_X Y, N\>  = \<A X, Y\> $,  and 

$H= \tr A = \sum_{i=1}^n h(e_i,e_i)$  for a local orthonormal frame $e_1, ..., e_n$  of the submanifold, where $\ona$ denotes the Levi-Civita connection on $\oM$. 

A Riemannian manifold with a density is a Riemannian manifold $(\oM, \og)$ where (without changing the metric) volumes  are measured with a weighted (smooth) function  $e^\psi: \oM \flecha \re$, in the following way: if $\Omega$ is a domain in $\oM$ and $M$ is a hypersurface with the induced metric $g$, the $(n+1)$-$\psi$-volume $V_\psi^{n+1}(\Omega)$ of $\Omega$ and  $n$-$\psi$-volume   $V^n_\psi(M)$ of $M$ are
\begin{equation*}
V_\psi^{n+1}(\Omega) = \int_\Omega e^\psi \ dv_\og,  \qquad V^n_\psi(M) = \int_{M} e^\psi \ dv_g,
\end{equation*}  
where $dv_\og$ and $dv_g$ are the $(n+1)$ and the $n$-volume elements induced by $\og$ and $g$ in the usual way on $\oM$ and $M$ respectively.   Obviously we have the corresponding volume  elements induced by the density
\begin{equation}\lb{dvp}
dv_\psi^{n+1} = e^\psi \ dv_\og , \qquad |dv^n_\psi| = e^\psi \ |dv_g| = e^\psi \ |\iota_N dv_\og| = |\iota_N dv^{n+1}_\psi|.
\end{equation}

Gromov (\cite{gro}) studied manifolds with densities as \lq\lq mm-spaces'', and mentioned the natural generalization of mean curvature in such spaces obtained as the gradient of the functional $n$-$\psi$-volume. According to \cite{gro}, \cite{mo} and \cite{rocabamo} it is denoted by $H_\psi$ and given (when $\ona\psi$ has sense) by 
\begin{align}\lb{Hpsi} H_\psi = H - \<\ona\psi, N\>.
\end{align}

We remark that this $H_\psi$ differs from the mean curvature $H'$ associated to the conformal metric $e^{2 \psi/n} \<\cdot, \cdot \>$ by a conformal factor, $H'= e^{-\psi} H_\psi$.

The geometry of manifolds with densities has received an increasing attention in the last years. In particular, in the context of extrinsic geometry, one can  see  the works \cite{canro,chmezh12a,chzhb,es, liwe,mamo,mo,mohoha,rocabamo} and references therein.

Hypersurfaces of $\re^{n+1}$ with $H_\psi=0$, called {\it $\psi$-minimal hypersurfaces (or $\psi$-minimal curves when $n=1$)}, appear  as self-similar solutions of the MCF when $e^\psi$ is a Gaussian, that is $\psi(x) =  \frac{C}{2} |x|^2$ (shrinkers if $C<0$ and expanders if $C>0$) or $\psi$ is linear (translating solitons). Other densities can appear in the study of MCF; for instance, other $\psi$-minimal curves in $\re^2$ appear when we look for shrinkers which are revolution hypersurfaces in $\re^3$ (see \cite{an92}). 

MCF is related with minimal surfaces and isoperimetric problems. These are mainly the kind of problems studied in the citations about manifolds with density quoted before. Then, it is  natural to consider,  in the setting of a manifold with density, a mean curvature flow governed by $H_\psi$ instead of $H$. We shall call this flow
\begin{equation}\label{gmcf}
\parcial{F}{t} =  \vec{H_\psi} = H_\psi\ N =  \(H - \<\ona\psi, N\>\)\ N,
\end{equation}
the {\it mean curvature flow with density $\psi$} (\DMCF for short).

In \cite{sm}, Smoczyk observed the equivalence between \DMCF of a hypersurface $M$ in $\oM$ and the MCF of $M\times \re$ in a warped product of $\oM$ and $\re$.
As a consequence, the \DMCF of a curve appears in a natural way in the study of the evolution of rotationally symmetric surfaces under MCF, like in \cite{sm94} and \cite{aag}. Variants of it (with some forcing term) appear in the study of the motion of hypersurfaces of revolution under volume preserving mean curvature flow (\cite{at97,at03,atka,cami09,cami12}).  

By a {\it radial density} on $\re^{n+1}$ we understand a density $e^\psi$ which is the composition of a function (that we will still  denote by the same letter) $\psi: [0,\infty[ \flecha \re$ and the distance to the origin $r(x)=|x|$. For radial densities $\psi$ that satisfy a special property, the flow of strongly starshaped hypersurfaces under \eqref{gmcf} was studied (with no explicit mention of the word \lq\lq density'') by Schnürer and Smoczyk in \cite{scsm}. Angenent (cf. \cite{an90} and \cite{an91}) and Oaks (cf. \cite{oa}) studied the  evolution of curves in a surface under a more general flow. The special case of the flow of convex surfaces in $\re^n$ with Gaussian and anti-Gaussian densities was studied by A. Borisenko and the first author in \cite{bomi6}.

Our aim is to start with a systematic study of the \DMCF. In this paper, after general considerations, we  study in some detail the \DMCF in a surface, specially, in the plane with a radial density. 

 A source of motivation has been the paper \cite{scsm}. There, the authors study the evolution of a strongly radial hypersurface in $\re^{n+1}-\{0\}$ with a radial density $\psi$ that satisfies some conditions. Our observation is that the most important properties of these conditions are i) the graph of $\psi'(r)$ cuts the graph of $-n/r$ only once and ii) before the crossing point, the graph of $\psi'$ is below that of $-n/r$. This seems to mean that the density with $\psi=- \ln r^n$ plays a special role among all the radial densities. This is confirmed by the results in \cite{bomi6}. Our first small contribution in this paper (section \ref{SVAD}) is the remark that the crossings between the graphs of $\psi'$ and $-n/r$ determine the spheres centered at $0$ which are $\psi$-minimal, and these are repulsors or attractors for the \DMCF depending on the relative position of the graphs in the neighborhood of the crossing points. This remark allows to determine the closed $\psi$-minimal hypersurfaces in some situations (Proposition \ref{psimin}).
 
 It is in the plane $\re^2$ with a radial density  where the above remarks (together with the results of Angenent and Oaks \cite{an90,an91,oa}) allow us to give a fairly complete description of the motion under \DMCF. In fact, in this context, we shall prove
 
 \begin{teor}\lb{pgen} Let 
$\re^2$ be the euclidean plane with  a radial density $\psi$ such that the graphs of  $\psi'$ and $-\fracc{1}{r}$ intersect transversally in  a discrete number of points $r_1 < r_2 < ... $,  let $\gamma_0$ be a simple closed curve which bounds a domain $\Omega_0$. Let $r_{max}$ and $r_{min}$ be, respectively, the maximal and minimal distance from $\gamma_0$ to the origin. Let us suppose that either the sequence of zeros $r_n$ goes to $\infty$, or the curve $\gamma_0$ is contained in the  disk centered at the origin whose radius is the biggest zero of  $\psi'+\fracc{1}{r}$, or  $\psi'+\fracc{1}{r}>0 $ after the biggest zero of  $\psi'+\fracc{1}{r}$. 

\noindent (a) If $\psi$ is smooth on $\re^2$, then:
\begin{itemize}
\item [a.i] If $r_{max}\le r_1$,  under \DMCF, $\gamma_0$ collapses to a point in finite time. In particular, if $\psi'+\fracc{1}{r}$ has no zero, every simple closed curve collapses to a point in finite time.
\item[a.ii]  If  $r_{2k-1} \le r \le r_{2k+1}$, $k\ge 1$, and $0\notin \Omega_0$,  under \DMCF, it collapses to a point in finite time.
\item [a.iii] If $r_{1} \le r_{min}$ and $0\in \Omega_0$,  the solution of the \DMCF with the initial condition $\gamma_0$ exists for $t\in[0,\infty[$ and, for every $m$, there is a sequence of times $t_n$, $t_n\to\infty$, such that the curves $\gamma(\cdot,t_n)$ converge, in the $C^m$ topology, to a $\psi$-minimal curve. Moreover:
\begin{itemize}
\item[a.iii.1] If $r_{2k-1} \le r \le r_{2k+1}$, $k\ge 1$,   the limit $\psi$-minimal curve is the circle of radius $r_{2k}$. This includes also the case $r_{2k+1}=\infty$, which occurs when $r_{2k}$ is the last zero of $\psi'+\fracc{1}{r}$. 
\end{itemize}
\end{itemize}

\noindent (b) If $\psi$  is smooth only on $\re^{2}-\{0\}$, $\lim_{t\to 0}\psi'(t) = -\infty$ and $\gamma_0$ is contained in $\re^2-\{0\}$:   
\begin{itemize}
\item [b.i] If $\psi'(r)>-1/r$ for $r<r_1$, the situation is the same as  in cases a.ii and a.iii.
\item [b.ii] If $\psi'(t)<-1/r$ for $r<r_1$, then: 
\begin{itemize}
\item[b.ii.1] If $\gamma_0$ is contained  inside the disk $r\le r_2$ and $0\notin \Omega_0$,  under \DMCF, it collapses to a point in finite time.
\item[b.ii.2]  If $\gamma_0$ is contained  inside the disk $r\le r_2$ and $0\in \Omega_0$,  the solution of the \DMCF with the initial condition $\gamma_0$ exists for $t\in[0,\infty[$, and, for every $m$, there is a sequence of times $t_n$, $t_n\to\infty$ such that the curves $\gamma(\cdot,t_n)$ converge, in the $C^m$ topology, to the circle of radius $r_1$.
\item [b.ii.3] The situations a.iii of the regular case are repeated but subtracting $1$ to the subindices of $r$.
\end{itemize}
\end{itemize}
\end{teor}

The following figure describes some of the possibilities given in this theorem. All the circles in the picture represent $\psi$-minimal circles. 
\begin{center}
\includegraphics[scale=0.4]{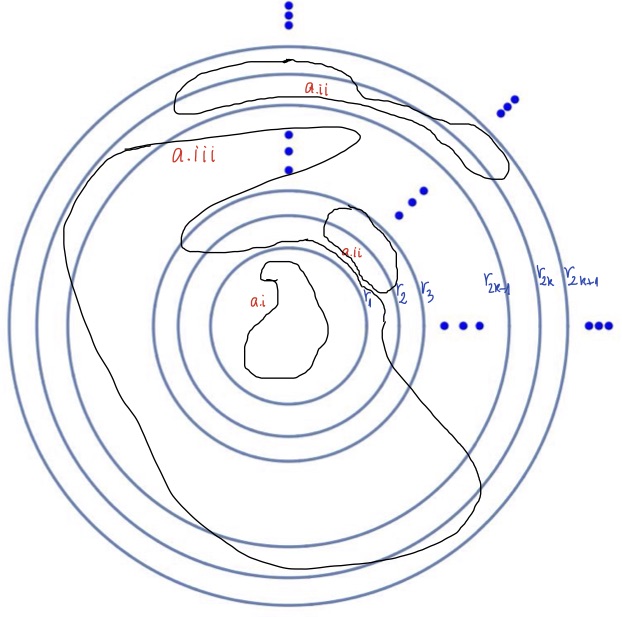}
\end{center}

The condition  \lq the graphs of  $\psi'$ and $-\fracc{1}{r}$ intersect transversally'' is not essential to give a description of the motions as we did, it only makes easier to describe the motion. After looking at the proof of the theorem, any reader can state easily other theorems where the condition is only that the number of intersections in any finite interval $[0,r_0]$ is finite, although the crossing could be only tangential and not transversal.

Let us observe that for the anti-Gaussian density, all the possibilities are included in case (ai), because $\psi'+\fracc{1}{r}$ has no zeros. The Gaussian density does not satisfy the condition \lq\lq $\psi'+\fracc{1}{r}>0 $ after the biggest zero of  $\psi'+\fracc{1}{r}$'', but we can still apply case (ai). We did not tray to go further in the study of this \DMCF with the methods of the proof of Theorem \ref{pgen} because Gaussian and anti-Gaussian densities are very special, and more complete results can be obtained by application of the method used in \cite{bomi6}. We shall do it in the appendix.

The situation studied by Schnürer and Smoczyk in \cite{scsm} is the one given by cases (b.ii) in the above theorem, with $r_2=\infty$. Here, cases (b.ii.1) and (b.ii.2) give all the possibilities of motion for a curve, then our results extend the result in \cite{scsm} for curves, where we do not need the extra hypothesis for the curve of being strongly starshaped.

Theorem \ref{pgen} confirms the special role played by the density $\psi= - \ln r^a$. The combinations of this density with Gaussian and anti-Gaussian,  $\psi(x)= \lambda  \fracc{|x|^2}{4} - a  \ln(|x|)$, with $\lambda\ne 0$, are also included in the situations described by Theorem \ref{pgen}. Among these combinations, there is one for which we can precise the  time of existence of the flow. This case, and the critical case $\psi(x)=- \ln(|x|)$, which does not fit in any of the situations described in Theorem \ref{pgen}, is the content of our second theorem.

\begin{teor}\lb{plnr} In the euclidean plane $\re^2$ with the density $\psi$, given any closed simple curve $\gamma_0$ contained in $\mathbb{R}^{2}-\lbrace0\rbrace$ that bounds a domain $\Omega_0$ of area $A$:
\begin{enumerate}\item  If $\psi(x)= - \ln(|x|)$ and $0\notin \Omega_0$, under \DMCF $\gamma_0$ evolves to a point in time $T=A/2\pi$.
\item  If $\psi(x)= - \ln(|x|)$ and $0\in \Omega_0$, the \DMCF with $\gamma_0$ as initial condition has solution for $t\in[0,\infty[$ and, for every $m$, there is a sequence of times $t_n$, $t_n\to\infty$, such that the curves $\gamma(\cdot,t_n)$ converge, in the $C^m$ topoloy, to a circle of radius $r=\sqrt{A/\pi}$.
\item If $\psi(x) = \lambda \frac{|x|^2}{4}-a \ln(|x|)$, with $\lambda >0$ and $a > 1$, and $0\notin \Omega_0$, 
 under \DMCF, $\gamma_0$ evolves to a point in time $T=\dfrac{1}{\lambda}ln\Big(1+\dfrac{\lambda A_{\overline{g}}(0)}{2\pi}\Big)$.
\end{enumerate}
\end{teor}

 For proving the above theorems we need some general facts remembered or studied in section \ref{sec2}. 
 
 Among these facts is important the variation of the area under \DMCF in dimension $2$. It is an unexpected fact for us that, although we are working with a flow driven by the gradient of the $\psi$-length, what matters for the study of it is the variation of the riemannian area, not the $\psi$-area. 
 
 Other basic tools are the results of Angenent and Oaks, and a study of the convergence to a $\psi$-minimal curve when the flow exists for all time in surfaces with general densities, which is the content of the following theorem, which has independent interest.
 It deals with curves evolving under \DMCF in a smooth Riemannian 2-manifold with density $(\oM^{2},g,e^{\psi})$ and Gauss curvature $K$ which, along the subset where the curve evolves, satisfies 
 \begin{align}
&\text{\rm i) }\  \begin{matrix}
|\ona^jK|\leq C_{j},\quad \vert \ona^j \psi\vert\leq P_{j}, \quad   \quad 0< E\leq e^{\psi} \le D \\
 \text{for some constants } C_{j}, P_{j}, E, D;\ j=0,1,2,...
 \end{matrix}  \lb{hipM1} \\
& \text{\rm ii) } \text{The isoperimetric profile $\mathcal I$ is a well defined continuous function which satisfies} \nn \\
& \qquad \text{ $\lim_{a\to a_0}\mathcal{I}(a) =0$ implies $a_0=0$}. \lb{hipM2} 
\end{align}
We remark that condition ii) is satisfied by every compact surface $\oM$ and many noncompact surfaces, among  that are those which cover a compact. In particular, $\re^2$ satisfies \eqref{hipM2}.

 \begin{teor}\lb{ctapg}
Let $(\oM^{2},\og,e^{\psi})$ be a $2$-Riemannian manifold with density satisfying \eqref{hipM2}. Let 
 $\gamma(\cdot,t)$ be a solution of the \DMCF \eqref{gmcf} with initial condition an embedded curve $\gamma_{0}:\mathbb{S}^{1}\flecha \oM^{2}$. If this solution exists for every $t\in[0,\infty[$, and $\gamma(\mathbb{S}^{1},t)$ is contained in a fixed compact domain $U$ where the conditions  \eqref{hipM1} are satisfied, then there is a reparametrization $\wga(\cdot,t)$ of $\gamma(\cdot,t)$ such that for every $m\in \ene$, there is a sequence  $\lbrace\wga(\cdot,t_{k})\rbrace_{k\in\mathbb{N}}$, $t_k\to\infty$,  which $C^m$-converges to a closed $\psi$-minimal curve of $\oM^2$.\end{teor}
 
 Let us remark that the hypothesis \lq\lq and $\gamma(\mathbb{S}^{1},t)$ is contained in a fixed compact domain $U$" is satisfied each time that we have a barrier for the \DMCF bounding a domain that contains $\gamma_0$, which happens in the situations described in theorems \ref{pgen} and \ref{plnr}.
 
The proof of Theorem \ref{ctapg} is the main technical part of the paper. The techniques used to prove this theorem are similar to those used for the MCF in \cite{ga1,ga2,gh,gr89} (see also \cite{chzh}).

The plan of the paper is as follows: In section \ref{sec2} we shall state the preliminaries about densities, with special emphasis on some variation formulae and the reasons why the \lq\lq critical" density $\psi(x)= - \ln|x|$ is important in the study of radial densities. The last include the classification of compact $\psi$-minimal hypersurfaces in some special circumstances.   This section finishes with the statement of the results of Angenent-Oaks on geometric flows that we shall need and under the form that we shall need. In section \ref{sec3} we shall give the proof of Theorem \ref{ctapg}, whereas theorems \ref{pgen} and \ref{plnr} will be proved in section \ref{sec4}. Finally, in the appendix, we shall do an independent and more complete study of the curve shortening flow for curves in the plane with the Gaussian and anti-Gaussian densities.

\section{Preliminaries}\lb{sec2}

 As a first important general property of \DMCF, we remark that  the same arguments used in MCF prove the following {\bf \lq\lq avoidance principle'' for \DMCF}: 
  
 {\it If $M$ and $\mathcal M$ are two smooth hipersurfaces, one of them compact, which do not intersect, the hypersurfaces obtained from them by \DMCF do not intersect. Moreover, embeddedness is preserved until the first singularity appears}.

\subsection{The $\psi$-divergence} \lb{sP1}
In a Riemannian manifold $\mathcal M$ with density $e^\psi$ we have, not only a special definition of the mean curvature, but also of the divergence and the laplacian. They are defined as the divergence and laplacian associated to the volume form $dv_\psi = e^\psi dv_g$. The definitions are:
\bec
(\dv_\psi X) \ dv_\psi = \mathcal L_X dv_\psi \quad \text{ and } \quad \Delta_\psi f = \dv_\psi\grad f.
\eec 
They produce the following computational formulae:
\bec
\dv_\psi X = \dv X + d\psi(X), \quad \Delta_\psi f = \Delta f + d\psi(\grad f).
\eec
We shall use the following consequence of the Stokes' Theorem  for this divergence:

\noindent {\it Given an oriented compact domain $\oO$ of a Riemannian manifold $\oM$ with smooth boundary $\partial \oO$ and a vector field $X$ on $\oO$, let $N$ be the unit vector normal to $\partial\oO$ pointing outward, one has 
\begin{align}
&\int_{\oO} f \oDelta_\psi f \ dv^{n+1}_\psi=  - \int_{\oO} |\ona f|^2 \ dv^{n+1}_\psi + \int_{\partial {\oO}}  f \<\ona f, N\> dv^n_\psi. \lb{divtp2}
\end{align}
}

In this paper we shall consider hypersurfaces $M$ in riemannian manifolds $\oM$ with a density $e^\psi$. We shall distinguish the laplacian, $\psi$-laplacian, gradient and covariant derivative in $\oM$ from those in $M$ by using an overline on the corresponding symbols when we are on $\oM$.

\subsection{Variation formulae}\lb{VF}

Let us consider a general motion driven by the vector field $f\ N$, where $f:M\times [0,\eps[ \flecha \re$,
\bec \parcial{F}{t} (x,t) = f(x,t) \ N(x,t). \lb{evo.Ff}
\eec

The corresponding evolution formulae for some geometric invariants associated to the immersions $F(\cdot , t)$ that we shall use in this paper are:

For the  riemannian $n$-volume  and $\psi$-$n$-volume elements $dv_{g_t}$ and $dv^n_{\psi}$, and the corresponding global volumes   $V^n_g(M_t)$ and $\psi$-$n$-volume $V^n_\psi(M_t)$,
 \bec\label{evol.dV}
\parcial{}{t} dv_{g_t} = - f\ H dv_{g_t}, \quad \text{ \  } \quad \parcial{}{t} dv^n_{\psi} = - f\ \Hp dv^n_{\psi}.
\eec
 \bec\label{evol.volF}
\deri{}{t} V^n_g(M_t) =  - \int_M f \ H \ dv_{g_t},\quad \text{ \  } \quad \deri{}{t} V^n_\psi(M_t) =  - \int_M f \ H_\psi \ dv^n_{\psi} .
\eec
Moreover
\begin{align}\lb{evol.onap}
\parcial{\<\ona \psi, N\>}{t} = \<\ona_{\parcial{F}{t}} \ona\psi, N\> + \<\ona \psi, \ona_{\parcial{F}{t}}N\>= f \<\ona_{N} \ona\psi, N\> - \<\ona \psi,  \nabla f\>.
\end{align}
\bec \lb{evolFH}
\parcial{H}{t} = \Delta f + f \Big(|A|^2+ \oRic(N,N)\Big). \eec
Joining \eqref{evol.onap} with  \eqref{evolFH}, we obtain
\begin{align} \lb{evolFHp}
\parcial{\Hp}{t} &= \Delta f + \<\ona\psi, \nabla f\>+ f \Big(|A|^2+ \oRic(N,N) -\<\ona_N\ona \psi, N\>\Big) \nn \\
&=\Delta_\psi f + f \Big(|A|^2+ \oRic_\psi(N,N)\Big),
\end{align}
where $ \oRic_\psi(N,N) :=  \oRic(N,N) -\<\ona_N\ona \psi, N\>$.

When  every $M_t = F_t(M)$ is embedded and the boundary of a compact domain $\Omega_t$ in $\oM$, we shall denote by $V^{n+1}_\og(\Omega_t)$ the $(n+1)$-volume  of $\Omega_t$. Its variation is given by
\begin{equation}\lb{igu_vol}
\deri{}{t} V^{n+1}_\og(\Omega_t) =  \int_{M} f\ \iota_{N} dv^{n+1}_\og = - \int_{M} f\ |dv^{n}_g|, \quad  \text{\rm if $N$ points inward}.
\end{equation}
In fact, we identify $M$ with $M_0=F_0(M)$ and consider a smooth (on $x$ an $t$) extension of $F_t$ to $\widetilde F_t:\Omega \equiv \Omega_0 \flecha \oM$
\begin{align*}
\deri{}{t} &V_\og(\Omega_t) = \deri{}{t} \int_\Omega  \widetilde F_t^*  dv_\og 
 = \int_{\Omega} \parcial{}{t} \(\widetilde F_t^* dv_\og \) = \int_\Omega   \mathcal L_{\parcial{\widetilde F_t}{t}} \widetilde F_t^*dv_\og \\
&=  \int_\Omega d \(\iota_{\parcial{\widetilde F_t}{t}} \widetilde F_t^*dv_\og \) = \int_M \iota_{\parcial{\widetilde F_t}{t}} \widetilde F_t^*dv_\og = \int_M   F_t^*\( \iota_{\parcial{F_t}{t}} \ dv_\og \),
\end{align*}
which is the claimed formula.
For the evolution of the distance to the origin, we have
\bec\label{evol.r}
 \parcial{r_t}{t} = \<\ona r,\parcial{F}{t} \>= f \<\ona r, N\>.
 \eec
 When $\oM=\re^{n+1}$, it is well known that 
 $$\Delta r = \fracc{n}{r}-\frac{|\nabla r|^2}{r} + H \< \ona r, N\> = \fracc{n}{r}-\frac{|\nabla r|^2}{r} + H_\psi \< \ona r, N\> + \<\ona \psi, N\> \< \ona r, N\>.$$
 When $f=H_\psi$, by substitution of this expression in \eqref{evol.r} one gets
 \bec\label{evol.re}
 \parcial{r_t}{t} = \Delta r -\fracc{n}{r}+\frac{|\nabla r|^2}{r} -\<\ona \psi, N\> \< \ona r, N\>.
 \eec

 \medskip

 \subsection{Second variation of the area with density and stability}\lb{SVAD}
From \eqref{evol.dV}, \eqref{evol.volF}  and \eqref{evolFHp} we can obtain the second variation for $V^n_\psi(M)$
\begin{align}\label{sec.dVpe}
 \deri{^2}{t^2} V^n_\psi(M_t) &=  - \int_M \deri{H_\psi}{t} \ f \ dv^{n}_{\psi} - \int_M H_\psi \ \deri{(f\ dv^n_{\psi})}{t} \nn \\
 &= - \int_M \(\Delta_\psi f + f \Big(|A|^2+ \oRic_\psi(N,N)\Big)\) f \ dv^n_\psi - \int_M H_\psi \ \deri{(f\ dv^n_{\psi})}{t}.
\end{align}
If $M_0$ is $\psi$-minimal, we obtain for its second variation (using \eqref{divtp2}):
\begin{align}
 \deri{^2}{t^2} V^n_\psi(M_0) 
 &=  \int_M \(|\nabla f|^2 - f^2 \ \( |A|^2+ \oRic_\psi(N,N)\) \) \ dv^n_\psi . \label{secvMinB}
\end{align}
In particular, when $\oM$ is the euclidean space $\re^{n+1}$, $\oRic_\psi = -\ona^2\psi$, and the formula of the second variation is 
\begin{align*}
 \deri{^2}{t^2} V^n_\psi(M_0) 
 & =  \int_M \(|\nabla f|^2 - f^2 \  \(|A|^2 - \ona^2\psi(N,N)\) \) \ dv^n_\psi . 
 \end{align*}
 
 Then, a $\psi$-minimal surface is strictly stable (attractor)  iff $ \(|A|^2 - \ona^2\psi(N,N)\) < 0$ and strictly inestable (repulsor) iff $\int_M  f^2 \ \(|A|^2 - \ona^2\psi(N,N)\) dv^n_\psi >\int_M|\ona f|^2\ dv^n_\psi$ for every variation $f$.
 
 We can consider the case of $\psi$-minimal hypersurfaces of the form $\psi(x)=$constant. For these hypersurfaces, the second fundamental form (when the chosen orientation is $N=\ona\psi/|\ona\psi|$) has the expression (for instance, see \cite{th} page 58):
 \bec
 h= - \frac1{|\ona\psi|} \ona^2\psi \text{ restricted to }TM.
 \eec
 Therefore, the hypersurface $\psi=$constant is $\psi$-minimal if and only if 
 \begin{align}
 &\fracc1{|\ona\psi|}\tr\ona^2\psi +\<\ona\psi,N\>=\frac1{|\ona\psi|}\ona^2\psi(N,N),\nn 
\end{align} 
that is (do not forget that  $N=\fracc{\ona\psi}{|\ona\psi|}$ on these hypersurfaces),
\begin{align}
 & \tr\ona^2\psi + |\ona\psi|^2=\ona^2\psi(N,N), \lb{pcm}
 \end{align}
 moreover, in general,  for these hypersurfaces we have
 \begin{align}
 |A|^2 = \fracc{|\ona^2\psi|^2}{|\ona\psi|^2} - \fracc{\ona^2\psi(N,N)^2}{|\ona\psi|^2} - 2 \sum_{i=1}^n \frac{(e_i|\ona\psi|)^2}{|\ona\psi|^2},
 \end{align}  
 because $\ona^2\psi(e_i,N)= e_i|\ona\psi|$.
 
So that, a $\psi$-minimal  $\psi=$constant hypersurface is strictly stable iff:
 \begin{align}\lb{stpct}
\fracc{|\ona^2\psi|^2 -\ona^2\psi(N,N)^2- 2 \sum_{i=1}^n (e_i|\ona\psi|)^2}{|\ona\psi|^2} - \ona^2\psi(N,N)<0.
 \end{align}
 
  \medskip
 When $\psi$ is  radial,
 
 $\ona \psi= \psi' \ona r$, $\ona^2\psi = \ona(\psi' \ona r) = \psi'' \ona r \otimes \ona r + \psi' \ona^2r$, $N= \fracc{\psi'}{|\psi'|}\ona r$.
 
In this case the conditon of stability \eqref{stpct} becomes 
\begin{align}\lb{stpctr}
\fracc{{\psi''}^2 + \psi'^2 \sum_{i,j=1}^n|\ona^2r(e_i,e_j)|^2  - {\psi''}^2}{{\psi'}^2} - \psi''<0 \quad \text{ i.e. }\quad  \psi'' > \sum_{i,j=1}^n|\ona^2r(e_i,e_j)|^2.
 \end{align}
 
 And the condition of $\psi$-minimality becomes:
\begin{align}
\fracc{\psi'' +  \psi' \sum_{i=1}^n \ona^2r(e_i,e_i)}{|\psi'|} + |\psi'| =\fracc{\psi''}{|\psi'|}, \text{ that is, } \psi' \sum_{i=1}^n \ona^2r(e_i,e_i) +  |\psi'|^2=0.
\end{align}

When we consider the less restrictive concept of stability under the \DMCF, with $\psi$ radial, it is easier to study the behavior of $\psi$-minimal spheres as attractors and repulsors, as the following observations show. 

Let us consider the {\bf evolution for the \DMCF of a sphere  $\mathbb{S}^{n}$ in $\re^{n+1}$} centered at the origin. If we choose $N=\ona\psi/|\ona\psi|$, we have  $\<\ona r, N\> =  \sgn(\psi')$ and $H_\psi =  - \sgn(\psi')(\frac{n}{r} +\psi')$, and formula \eqref{evol.r} becomes 
 \begin{align}
 \parcial{r}{t} &=  - \sgn(\psi')(\frac{n}{r} +\psi') \sgn(\psi') = -(\frac{n}{r} +\psi'). \lb{drdtS}
 \end{align}
 Then, 
 \begin{prop}\lb{cilindros} In $\re^{n+1}$ with a radial density, a sphere 
 $\mathbb{S}^{n}$ of radius $r$ and centered at the origin is $\psi$-minimal if and only if $\psi'=-\fracc{n}{r}$.

Under the \DMCF, the radius $r$ of a sphere $\mathbb{S}^{n}$ centered at the origin

 - increases when  $- (\fracc{n}{r} +\psi') >0$,  that is,   when  $\psi'<-\fracc{n}{r}$,
 
- decreases when $-(\fracc{n}{r} +\psi') <0$,  that is, when   $\psi'> -\fracc{n}{r}$ . 
 \end{prop}
 In the following pictures, the curve $(r, -n/r)$ in marked with black ink. They show two examples for the function $\psi'(r)$  (in red), one for a density with singularities and the other for a regular density. For each one of these examples, the spheres of radius in the verticals marked with A are attractors (other neighbor spheres are attracted by them) and those with radius in the verticals marked with R are repulsors. Those marked with $A$ or $R$ are $\psi$-minimal and, therefore, stationary. The arrows indicate the direction in which a  sphere of radius equal to the coordinate \lq\lq$r$'' of the basis of the arrow  moves.
 
 \begin{center}
    
\includegraphics[scale=0.4]{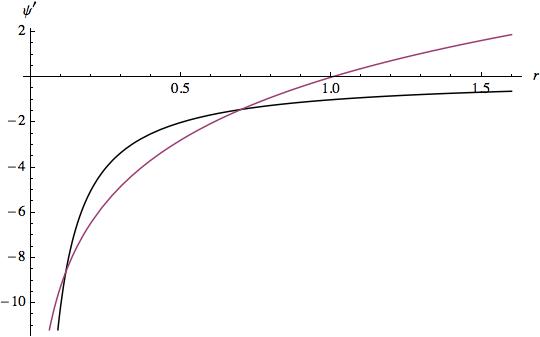}
\includegraphics[scale=0.4]{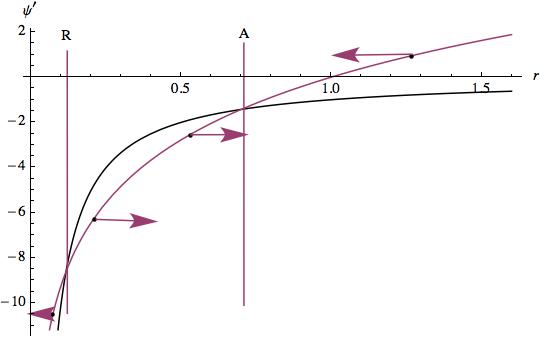}

Picture 1

\includegraphics[scale=0.4]{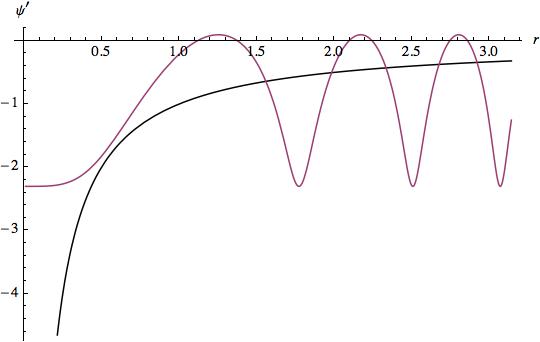}
\includegraphics[scale=0.4]{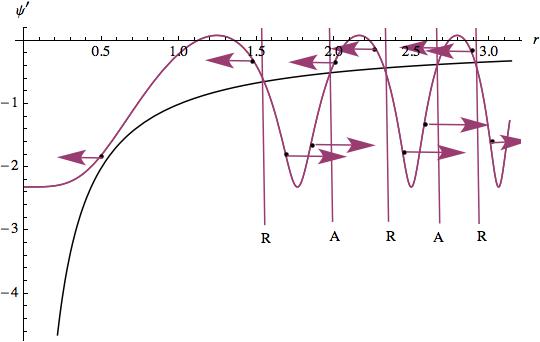}

Picture 2

 \end{center}

As a consequence of this behavior and the avoidance principle for \DMCF, we obtain:

\begin{prop}\lb{psimin} 
Let $\re^{n+1}$ be the euclidean space with a  radial density $\psi$.   Let us denote  by $r_{max}$ (respectively $r_{min}$) the maximal (respectively minimal) distance from any point of a closed hypersurface $\Sigma$ to the origin. Let us suppose that 

(a) $\psi$ is smooth, and 

(b) the graphs of  $\psi'$ and$-\fracc{n}{r}$ intersect transversally in  a discrete family of points $r_1 < r_2 < ... $. Then
\begin{itemize}
\item [i]  A hypersurface $\Sigma$ satisfying  $r_{2k} \le r_{max}\le r_{2k+1}$,  is $\psi$-minimal if and only if it is one of the spheres $r=r_{2k}$ or $r= r_{2k+1}$. In particular, there is no closed $\psi$-minimal hypersurface  inside the sphere $r=r_1$.
\item [ii] The unique closed $\psi$-minimal  hypersurfaces $\Sigma$ with $r_{2k+1} \le r_{min}\le r_{2k+2}$ are the spheres $r=r_{2k+1}$ and $r= r_{2k+2}$ .
\item [iii] If the number of zeros of $\psi'+\fracc{n}{r}$ is even and $r_k$ is the last one, the unique closed $\psi$-minimal  hypersurface $\Sigma$ with $r_{max}\ge r_k$ is the sphere $r=r_k$.
\item [iv] If the number of zeros of $\psi'+\fracc{n}{r}$ is odd and $r_k$ is the last one, the unique closed $\psi$-minimal  hypersurface $\Sigma$ with $r_{min}\ge r_k$ is the sphere $r=r_k$.
\item[v] When $n=1$, a $\psi$-minimal simple closed curve needs to be starshaped respect to the origin.
\end{itemize}

If we change the hypothesis (a) by 
(a')  $\psi$  is smooth only on $\re^{2}-\{0\}$, $\lim_{t\to 0}\psi'(t) = -\infty$ and $\Sigma$ is contained in $\re^{n+1}-\{0\}$:
\begin{itemize}
\item[i'] If $\psi'(r)>-\fracc{n}{r}$ for $r<r_1$,  the same situations than under the hypothesis (a) are repeated.
\item[ii'] If $\psi'(r)<-\fracc{n}{r}$ for $r<r_1$, we have again the same situations than in case $\psi$ smooth, but interchanging the roles of $r_{max}$ and $r_{min}$.
\item[iii'] Case v in the smooth situation holds with no change.
\end{itemize}
\end{prop}
\begin{demo}
In each one of cases (i) to (iv), (i') and (ii'), if $\Sigma$ is $\psi$-minimal and is not one of the declared spheres, there is a sphere disjoint with $\Sigma$ which evolves under \DMCF until it touches $\Sigma$, which is in contradiction with the avoidance principle.

Let us considere the case (v). In a radial density, the lines through the origin are $\psi$-minimal. On the other hand, $\Sigma$ starshaped respect to the origin means that $\<F,N\> \ne 0$ at every point $F$ of $\Sigma$. Then, if $\Sigma$ is not starshaped, there is a line through the origin tangent to $\Sigma$  but, by the maximum principle, this is impossible. 
\end{demo}

For the critical case $\psi(x) = a - n \ln(|x|)$, it is possible to give a complete classification of the closed $\psi$-minimal hypersurfaces:

\begin{prop}\lb{critmin}
In $\re^{n+1}$ with the critical density $\psi(x) = a - n \ln(|x|)$, the unique closed $\psi$-minimal hypersurfaces are the spheres centered at the origin
\end{prop}
\begin{demo}
We know that all the spheres centered at the origin are $\psi$-minimal. If $\Sigma$ is another closed $\psi$-minimal surface, it must be tangent to some of these spheres at some point but, because they are two tangent $\psi$-minimal surfaces, by the maximum principle, they must coincide. 
\end{demo}

For the Gaussian case $\psi(x)=-1/2\ \mu^2 r^2(x)$ in the plane, it is known (cf. \cite{AL}) that the unique  $\psi$-minimal closed simple curve is the circle of radius $r=1/\mu$. For higher dimensions, it is an active research field to find all the Gaussian minimal surfaces, which coincide with the shrinkers of the MCF (see, for instance, \cite{br,moll}).

\subsection{Angenent-Oaks' results}

Given a riemannian surface $\oM$, let $S^1\oM$ denote the unit tangent bundle of $\oM$. Angenet (\cite{an90,an91}) studied the solution $F(\cdot,t):\mathbb{S}^{1}\flecha \oM$ of the flow 
\begin{equation}
\<\parcial{F}{t},N\> = V(\ft,H), \quad   \lb{Vmcf}
\end{equation}
 for $V:S^1\oM\times\re \flecha \re$ satisfying some conditions called V1), V2), V3), V4) and V5*). The case of \DMCF on $\oM$ corresponds to the special choice of $V(\ft_x, \ell)= \ell - \<\ona \psi(x),\ft_x\>$. For this $V$, condition V2) is automatically satisfied, and the others can be written in the following form: There are positive numbers $\mu, \hat\mu, \nu$ such that

\begin{itemize}
\item[V1)] $\ona\psi$ is a locally Lipschitz continuous function,
\item[V3)]  $|\ona\psi|\le \mu$ for almost all $x \in \oM$,
\item[V4)] $|\ona^2 \psi(\ft,\cdot)|^2+ \<\ona\psi, J\ft\>^2 \le \hat\mu^2$ for every $\ft\in S^1\oM$,
\item[V5*)] $|\ona^2 \psi(\ft,\cdot)|+|\ell|\ |\<\ona\psi, J\ft\>| \le \nu (1+|\ell|^2)$ for almost all $(\ft,\ell)\in S^1\oM\times \re$,
\end{itemize}
where $Jt\in S^1\oM$ is orthogonal to $t$.
Moreover this function $V(\ft_x, \ell)= \ell - \<\ona \psi(x),\ft_x\>$ also satisfies the following property used in \cite{oa}

V6) $V(-\ft,-H) = -H-\<\ona \psi, \ft\> = - V(\ft,H)$.

Let us remark that all these properties V1) to V6) are satisfied if $\psi$, $\ona\psi$ and $\ona^2\psi$ are smooth and bounded on $\oM$.

The main result of Angenent and Oaks that we shall need, particularized to the \DMCF on $\oM$ with a density $e^\psi$ satisfying that $\psi$, $\ona\psi$ and $\ona^2\psi$ are smooth and bounded is 

\begin{teor}{{\rm\cite{an90,an91,oa}}}\lb{ThAO}
 Let $F_0: \mathbb{S}^{1} \flecha \oM$ be a simple $C^2$ curve. Then the solution
$F(u,t)$ to \eqref{gmcf} such that $F(u,0) = F_0(u)$ either shrinks to a point on $\oM$ in finite
time or exists for infinite time.
\end{teor}
 This was stated as a corollary in \cite{oa}.
 
 \section{Convergence to a $\psi$-minimal curve in $\oM^2$.}\lb{sec3}
 
 In this section we shall prove Theorem \ref{ctapg}.
Since we are dealing with curves, we shall use the letter $\gamma:\mathbb{S}^{1} \flecha \oM$ instead of $F$ to denote the immersions from $\mathbb{S}^{1}$ into $\oM$, $k$ and $k_\psi$ to denote the corresponding curvatures $H$ and $H_\psi$. Moreover, we shall use $ds$, $L$ and $A$ instead of $dv^1$, $V^1$ and $V^2$. Because $n=2$, the Ricci curvature $\oRic = K\ \og$, and $\oRic_\psi = K\ \og - \ona^2\psi$. 

For the proof of Theorem \ref{ctapg} we shall require some preparatory lemmas.

\begin{lema}\lb{lema6}
  Under the hypothesis of Theorem \ref{ctapg}, the conditions 
\begin{align}
 &\text{i) there is a constant $c>0$ such that } L(t)\geq c \text{ for every } t\in[0,\infty[,\label{ProSol2} \\
 &\text{ii) the riemannian area enclosed by the evolving curves is bounded from above,}
\lb{ProSol4}
 \end{align}
  are satisfied.
 \end{lema}
\begin{prue} To prove claim i) we shall prove that there is a constant $c'$ such that \linebreak $\liminf_{t\rightarrow \infty}L(t)\geq c'$. This is a consequence of the following fact that we shall prove. \lq \lq Let $[0,T_{max}[$ be the maximal interval of $t$ for which the solution $\gamma(\cdot,t)$ of \eqref{gmcf} on $\oM^2$ is well defined. If $\ds 
\liminf_{t\rightarrow T_{max}}L(t)=0$, then  $T_{max}<\infty$." To prove it, first we observe that the evolution equation \eqref{evol.volF} implies that $L_\psi$ decreases with time. This property, together with the hyphotesis  $0<E\leq e^{\psi}\leq D$, gives the following chain of implications
\begin{align}
\liminf_{t\rightarrow T_{max}} L(t)=0 \Rightarrow \liminf_{t\rightarrow T_{max}} L_{\psi}(t)=0\Rightarrow \lim_{t\rightarrow T_{max}} L_{\psi} (t)=0 \Rightarrow \lim_{t\rightarrow T_{max}} L(t)=0.
\end{align}
 From the hypothesis \eqref{hipM2} it follows that $\lim_{t\rightarrow T_{max}} L(t)=0$ implies $\lim_{t\rightarrow T_{max}} A(t)=0$. Then, there is a $t_{0}\in[0,T_{max}[$ such that
\begin{align}\lb{Aotb}
{\max_{U}\big( |K|+|\overline{\Delta}\psi|\big)}\ A(\Omega_{t})\leq \pi \text{ \ for every } t\geq t_{0}.
\end{align}
From \eqref{igu_vol} it follows that 
\begin{align}
\dfrac{d}{dt}A&=-\int_{\mathbb{S}^{1}}k_{\psi}ds_{t}=-\int_{\mathbb{S}^{1}}k ds_{t}+\int_{\mathbb{S}^{1}}\<\overline{\nabla}\psi,N\>ds_{t} \nn
\\
&=-2\pi+\int_{\Omega}Kda-\int_{\Omega}\overline{\Delta}\psi da 
=-2\pi+\int_{\Omega}\big( K-\overline{\Delta}\psi\big)da \nn
\\
&\leq-2\pi+\max_{U}\big( |K|+|\overline{\Delta}\psi|\big)A(\Omega), \lb{datK}
\end{align}
From \eqref{Aotb} and \eqref{datK} we conclude that
$\ds\dfrac{d}{dt}A\leq -\pi$ for every $t\geq t_{0}$, which implies that $T_{max}$ is finite. 

On the other hand, the condition \eqref{ProSol4} is a direct consequence  of the hypotheses that $\gamma(\mathbb{S}^{1},t) \subset U$ compact for every $t\in [0,T_{max}[$.
\end{prue}

\begin{lema}\lb{le8}
Let $(\oM,g,e^{\psi})$ be a  surface with density satisfying \eqref{hipM1}, and let $\gamma$ be an embedded curve  evolving under \DMCF.  If $\gamma(\cdot,0)$ satisfies $k_\psi \ge C = P_1 + \sqrt{C_0+ P_2}$, then
$
k_{\psi}(\cdot,t)\geq C \:\text{ for every } t\in[0,T_{max}).
$
\end{lema}
\begin{proof}
From \eqref{evolFHp}, the variation formula of  $k_{\psi}$ is:
\begin{align*}
\dfrac{d}{dt}k_{\psi}=\Delta k_{\psi}+\<\overline{\nabla}\psi,\nabla k_{\psi}\>+k_{\psi}(k^{2}+K-\ona^2\psi(N,N)),
\end{align*}
but, because of the hypotheses \eqref{hipM1}, if $k_{\psi} \ge P_1 + \sqrt{C_0+ P_2}$, then  $k^{2}+K-\ona^2\psi(N,N)\geq (k_\psi+\<\ona \psi,N\>)^2-C_0- P_2  \ge (k_\psi-P_1)^2-C_0- P_2 = 0$; and the lemma follows from the maximum principle.
\end{proof}
\begin{lema} If $s_t$ denotes the arc-length parameter of $\gamma(\cdot,t)$, one has
\begin{equation}\lb{corst}
[\partial_{t},\partial_{s_{t}}]=k_{\psi}k\partial_{s_{t}}.
\end{equation}
\end{lema}
\begin{proof} If we parametrize the curve $\gamma(\cdot,t):\mathbb{S}^{1} \flecha \oM$ by the angle $\alpha$ of $\mathbb{S}^{1}$, one has
$
ds_{t}=ds_{t}(\partial_{\alpha})d\alpha
$,
then
$\ds
\partial_{s_{t}}=\dfrac{1}{ds_{t}(\partial_{\alpha})}\partial_{\alpha}
$
and, using formula \eqref{evol.dV}
\begin{align*}
[\partial_{t},\partial_{s_{t}}]&
=\partial_{t}\Bigg(\dfrac{1}{ds_{t}(\partial_{\alpha})}\Bigg)\partial_{\alpha}
=\dfrac{k_{\psi}k}{ds_{t}(\partial_{\alpha})^{2}}ds_{t}(\partial_{\alpha})\partial_{\alpha}
=\dfrac{k_{\psi}k}{ds_{t}(\partial_{\alpha})}\partial_{\alpha}=k_{\psi}k\partial_{s_{t}}. \nn
\end{align*}
\end{proof}

{\bf Proof of Theorem \ref{ctapg}}. We shall do it in some steps:

{\bf Step 1.} {\it Let  $\gamma:\mathbb{S}^{1}\times[0,\infty)\rightarrow \oM^{2}$ be a solution of the \DMCF satisfying (\ref{ProSol2}) and (\ref{ProSol4}) moving in a surface with density $(\oM^{2},g,e^{\psi})$  satisfaying \eqref{hipM1}.
Then, there is a constant $C>0$ such that 
$
\text{ for every } \quad t\in [0,\infty)  \quad \inf_{\mathbb{S}^{1}}\vert k_{\psi}(\cdot,t)\vert\leq C.
$
}

\begin{prue}
Let us suppose that, on the contrary, 
\begin{align*}
\text{for every }\:C>0,\:\text{ there is a }\: t_{0}=t_{0}\big(C\big)\:\text{ such that }\: \inf_{\mathbb{S}^{1}}\vert k_{\psi}(\cdot,t_{0})\vert > C.
\end{align*}
Let us choose an arbitrary $C>0$.
There are two possibilities
\begin{align}
k_{\psi}(s,t_{0})&=\vert k_{\psi}(s,t_{0})\vert>C>0,\:\forall s\in\mathbb{S}^{1}, \lb{posib1}
\\
&\:or \nn
\\
k_{\psi}(s,t_{0})&=-\vert k_{\psi}(s,t_{0})\vert< - C<0,\:\forall s\in\mathbb{S}^{1}. \lb{posib2}
\end{align}

 In case \eqref{posib1}, let us select $C = P_1+\sqrt{C_0+P_2}+1$. By Lemma \ref{le8},  
$
k_{\psi}(s,t)\geq C,\:\text{ for every } (s,t)\in\mathbb{S}^{1}\times[t_{0},\infty).
$
Therefore,
\begin{equation*}
\deri{}{t}L_{\psi}=-\int_{\mathbb{S}^{1}}k_{\psi}^{2}e^{\psi}ds\leq -C^2\ L_{\psi}\leq -C^2\ E\ L\leq -C^2\ E\ c<0,
\end{equation*}
where we have applied \eqref{hipM1} in the second inequality and \eqref{ProSol2} in the third one. This implies  $T_{max}<\infty$ in contradiction with the hypothesis that the solution $\gamma(\cdot,t)$ is well defined for every $t\in[0,\infty)$.

In case \eqref{posib2}, this inequality implies  $ k\leq P_1 -C$ for every $s\in\mathbb{S}^{1}$.
Now we choose  $C=P_1+n$, $n$ big enough, then  $k\le -n$. If we take into account also \eqref{ProSol2} and Gauss-Bonnet Theorem, we obtain 
\begin{align*}
-n\ c \geq-n\ L(t_{0})&\geq\int_{\mathbb{S}^{1}}k ds_{t_{0}}=2\pi-\int_{\Omega_{t_{0}}}Kda
\geq 2\pi-C_0 A(\Omega_{t_{0}}),
\end{align*}
which is a contradiction with \eqref{ProSol4} if we take $n$ big enough.

This finishes the proof of step 1.
\end{prue}

{\bf Step 2} {\it Let  $\gamma:\mathbb{S}^{1}\times[0,\infty)\rightarrow \oM^{2}$ be a solution of the \DMCF satisfying (\ref{ProSol2}) and (\ref{ProSol4}) moving in a surface with density $(\oM^{2},g,e^{\psi})$  satisfaying \eqref{hipM1}. 
If this solution exists for every $t\in[0,\infty[$, then 
\begin{equation}\label{P_kpsicef}
\lim_{t\rightarrow\infty}\int_{\gamma(\cdot,t)}k_{\psi}^{2}(s,t)ds_{\psi}=0.
\end{equation} 
}

\begin{prue}  For it, we need to compute the variation of the $L_\psi^2$-norm of $k_\psi$. We shall use  \eqref{hipM1}, \eqref{divtp2},  \eqref{evolFHp}, $\oRic_\psi = K \og - \ona^2\psi$ and $k^{2}= \(k_\psi + \<\ona \psi,N\>\)^2\leq \wkp + 2 |k_\psi| P_1 + P_1^2 \le (1+ 2 P_1)\wkp\ + P_1^2 + 2P_1$ to obtain
\begin{align}
\dfrac{d}{dt}\int_{\gamma(\cdot,t)}k_{\psi}^{2}ds_{\psi}&=\int_{\gamma(\cdot,t)}\Big[2k_{\psi}\partial_{t}(k_{\psi})-k_{\psi}^{4}\Big]ds_{\psi}
\nn\\
&=\int_{\gamma(\cdot,t)}\Big[2k_{\psi}\Big(\Delta_{\psi}k_{\psi}+k_{\psi}(k^{2}+\overline{Ric}_{\psi}(N,N))\Big)-k_{\psi}^{4}\Big]ds_{\psi}
\nn\\
&=\int_{\gamma(\cdot,t)}\Big[2k_{\psi}^{2}k^{2}-k_{\psi}^{4}-2(\partial_{s}k_{\psi})^{2}+2k_{\psi}^{2}\overline{Ric}_{\psi}(N,N)\Big]ds_{\psi} \nn \\
& \le a \int_{\gamma(\cdot,t)} k_\psi^{4} ds_{\psi}  +   b \int_{\gamma(\cdot,t)} \wkp ds_{\psi}-2\int_{\gamma(\cdot,t)}(\partial_{s}k_{\psi})^{2}ds_{\psi} , \lb {dtkp1}
\end{align}
for some positive constants $a$ and $b$ independent of $t$. 

In order to bound the first addend of the last inequality, let us observe that there are constants $a_1$, $b_1>0$ independent of $t$ such that
\begin{equation}\lb{bmaxkp}
\max_{\gamma(\cdot,t)}k_{\psi}^{2}\leq a_1 + b_1\int_{\gamma(\cdot,t)}(\partial_sk_{\psi})^{2}ds_{\psi},\:\forall t\in[0,\infty).
\end{equation}
In fact,
from Step 1, there is a $C>0$ such that
\begin{center} $\inf_{\mathbb{S}^{1}}\vert k_{\psi}(\cdot,t)\vert\leq C,\:\text{ for every } t\in [0,\infty)$.
\end{center}

Given any fixed $t\in[0,\infty)$, let $s_{0}\in\mathbb{S}^{1}$, $s_{0}=s_{0}(t)$, satisfying $\vert k_{\psi}(s_{0},t)\vert\leq C$. By integration,
\begin{align*}
\vert k_{\psi}(s_{1})\vert & \leq
\vert k_{\psi}(s_{0})\vert+\Big\vert\int_{s_{0}}^{s_{1}}\partial_{s}k_{\psi}ds_{t}\Big\vert\leq
C+\Big\vert\int_{s_{0}}^{s_{1}}\partial_{s}k_{\psi}ds_{t}\Big\vert
\nn\\
&\qquad \le \sqrt{2}\sqrt{C^2 + \(\int_{s_0}^{s_1}(\partial_{s}k_{\psi})^{2}ds_{t}\) \(\int_{s_0}^{s_1}ds_{t}\)}\leq
\sqrt{2}\sqrt{C^2+\dfrac{1}{E^2}L_{\psi}(t)\int_{\gamma(\cdot,t)}(\partial_{s}k_{\psi})^{2}e^{\psi}ds_{t}},
\end{align*}
where we have used Hölder in the third inequality and \eqref{hipM1} in the fourth one. Since $L_{\psi}(t)$ is decreasing, we obtain from the above expression that
\begin{equation*}
\max_{\gamma(\cdot,t)}k_{\psi}^{2}\leq 2  C^{2}+\frac{2 L_\psi(0)}{E^2}\int_{\gamma(\cdot,t)}(\partial_{s}k_{\psi})^{2}ds_{\psi},
\end{equation*}
which proves \eqref{bmaxkp}.

Plugging \eqref{bmaxkp} into \eqref{dtkp1} we obtain that there are positive constants $a_2, b_2$ independent of $t$ satisfying
\begin{align}
\dfrac{d}{dt}\int_{\gamma(\cdot,t)} & k_{\psi}^{2}ds_{\psi} \nn \\
&\leq 
-2\int_{\gamma(\cdot,t)}(\partial_{s}k_{\psi})^{2}ds_{\psi}+ a_2 \int_{\gamma(\cdot,t)}k_{\psi}^{2}ds_{\psi}+ b_2 \int_{\gamma(\cdot,t)}k_{\psi}^{2}ds_{\psi}
\int_{\gamma(\cdot,t)}(\partial_{s}k_{\psi})^{2}ds_{\psi}, \lb{dtab}
\end{align}
for every $t\in[0,\infty)$.

On the other hand, given an arbitrary $\eps$, with $0<\eps<\fracc{1}{b_2}$,  there is $t_{1}=t_{1}(\eps)$ such that 
\begin{equation} \lb{epseps}
 \int_{\gamma(\cdot,t_1)}k_{\psi}^{2}ds_{\psi} <\frac{\eps}{2} \quad \text{ and } \quad \int_t^\infty \int_{\gamma(\cdot,t)}k_{\psi}^{2}ds_{\psi}\ dt <\eps^2 \text{ for } t\in[t_1(\eps),\infty[.
 \end{equation}
 In fact, from the variation formula \eqref{evol.volF} we know that $L_\psi$ is decreasing with time, then 
\begin{align}
L_\psi(0) \ge - \lim_{t\to \infty} \int_0^t \deri{L_\psi}{t} dt = \lim_{t\to \infty} \int_0^t \int_{\gamma(\cdot,t)}k_{\psi}^{2}ds_{\psi}\ dt, \nn
\end{align}
and this implies \eqref{epseps}.

From inequality \eqref{dtab} we obtain that, for any $t'$ in which the inequality
 \begin{equation}\lb{condi}
 \int_{\gamma(\cdot,t')}k_{\psi}^{2}ds_{\psi}\leq 1/b_2
 \end{equation}
 is satisfied, one has
\begin{align}
\dfrac{d}{dt}\int_{\gamma(\cdot,t')}k_{\psi}^{2}ds_{\psi}&\leq
 -\int_{\gamma(\cdot,t')}(\partial_{s}k_{\psi})^{2}ds_{\psi}+ a_2 \int_{\gamma(\cdot,t')}k_{\psi}^{2}ds_{\psi}. \label{E_clavekpsicero}
\end{align}

Let us take $t_1=t_1(\eps)$ of inequality \eqref{epseps}. 
We claim that 
 $\int_{\gamma(\cdot,t)}k_{\psi}^{2}ds_{\psi}<\eps,\:\text{for every } t\geq t_{1}$.

In fact. If not, let  $t_{2}>t_{1}$ be the first $t$ for which
$
\int_{\gamma(\cdot,t_{2})}k_{\psi}^{2}ds_{\psi}=\eps.
$
So \eqref{condi} is satisfied for every $t'\in[t_1,t_2]$, and we can also use  \eqref{E_clavekpsicero} on the interval $[t_1,t_2]$ to obtain
\begin{align*}
\dfrac{\eps}{2}&< \int_{\gamma(\cdot,t_{2})}k_{\psi}^{2}ds_{\psi}-\int_{\gamma(\cdot,t_{1})}k_{\psi}^{2}ds_{\psi}=
\int_{t_{1}}^{t_{2}}\Big(\dfrac{d}{dt}\int_{\gamma(\cdot,t)}k_{\psi}^{2}ds_{\psi}\Big) dt 
\nn\\
& \leq\int_{t_{1}}^{t_{2}} a_2 \int_{\gamma(\cdot,t)}k_{\psi}^{2}ds_{\psi}\leq a_2\int_{t_{1}}^{\infty}\int_{\gamma(\cdot,t)}k_{\psi}^{2}ds_{\psi} dt < a_2 \eps^{2}.
\end{align*}
This implies
$\dfrac{1}{2 a_2}<\eps$,
but $\eps<\dfrac{1}{b_2}$ was arbitrary, which gives a contradiction, then our claim is true and we conclude \eqref{P_kpsicef}, which finishes the proof of Step 2.
\end{prue}


{\bf Step 3.}
\ {\it Let  $\gamma:\mathbb{S}^{1}\times[0,\infty)\rightarrow \oM^{2}$ be a solution of the \DMCF satisfying (\ref{ProSol2}) and (\ref{ProSol4}) moving in a surface with density $(\oM^{2},g,e^{\psi})$  satisfaying \eqref{hipM1}. 
If this solution exists for every $t\in[0,\infty[$, then 
\begin{equation}\lb{acotg}
\lim_{t\rightarrow\infty}\int_{\gamma(\cdot,t)}(\partial_{s}^nk_{\psi})^{2}ds_{\psi}=0 \ \text{ for every natural number }\ n.
\end{equation}
}
We shall prove it by induction. When $n=0$, \eqref{acotg} is just \eqref{P_kpsicef}. Let us suppose that \eqref{acotg} is true for the derivatives of $k_\psi$ until order $n-1$. We shall see that it is also true for the derivative of order $n$. As in the previous step, we start by computing the derivative respect to $t$.

\begin{align}
\dfrac{d}{dt}\int_{\gamma(\cdot,t)}(\partial_{s}^n k_{\psi})^{2}ds_{\psi}&=
\int_{\gamma(\cdot,t)}\Big[2 (\partial_{s}^n k_{\psi})\ \partial_{t}\partial_{s}^n (k_{\psi})-(\partial_{s}^n k_{\psi})^{2}k_{\psi}^{2}\Big]ds_{\psi}.  \lb{dtidns2a}
\end{align}
Using the commutation rule \eqref{corst}, we obtain 
\begin{align}
\partial_{t}\partial_{s}^n (k_{\psi}) = \partial_s^n \partial_t k_\psi + \sum_{i=0}^{n-1} \partial_s^i (k_\psi k \partial_s^{n-i}k_\psi) \nn
\end{align}
and, plugging here the variation formula \eqref{evolFHp} of $k_\psi$
\begin{align}\lb{dtdns2}
\partial_{t}(\partial_{s}^n k_{\psi}) =  \partial_s^n \(\Delta_\psi k_\psi + k_\psi \ (k^2 - K + \ona^2\psi(N,N)) \) + \sum_{i=0}^{n-1} \partial_s^i (k_\psi k \partial_s^{n-i}k_\psi).
\end{align}
Let us compute each one of the terms in \eqref{dtdns2}
\begin{align}
\partial_s^n (\Delta_\psi k_\psi) &= \partial_s^{n+2} k_\psi + \partial_s^n\<\nabla\psi,\nabla k_\psi\> \nn \\
&= \Delta (\partial_s^n k_\psi) + \partial_s \psi \ \partial_s^{n+1}k_\psi + \sum_{i=1}^n \binom{n}{i} \partial_s^{1+i} \psi \ \partial_s^{n+1-i}k_\psi \nn \\
& = \Delta_\psi (\partial_s^n k_\psi)+ n \ \partial_s^{2} \psi\  \partial_s^{n} k_\psi  +\sum_{i=2}^n \binom{n}{i} \partial_s^{1+i} \psi\ \partial_s^{n+1-i}k_\psi, \lb{ter1}
\end{align}

\begin{align}
\partial_s^n \( k_\psi \ (k^2 - \oRic_\psi(N,N)) \)&= (\partial_s^n k_\psi) (k^2 - K + \ona^2\psi(N,N))  \nn \\
& \ + \sum_{i=1}^{n} \binom{n}{i}\partial_s^i (k^2 - K + \ona^2\psi(N,N)) \partial_s^{n-i} k_\psi ,\lb{ter2}
\end{align}

\begin{align}
 \sum_{i=0}^{n-1} \partial_s^i (k_\psi k \partial_s^{n-i}k_\psi) &= k_\psi k \partial_s^nk_\psi + \sum_{i=1}^{n-1} \partial_s^i (k_\psi k \partial_s^{n-i}k_\psi) .\lb{ter3}
\end{align}
Since, from now on, the notation with indices is becoming complicated, we shall use the following multi-index notation. Capitals will denote multi-indices. For us, all the entries $j_k$ of a multi-index  $J=(j_1, ...., j_q)$ will be ordered $j_1 \ge j_2 \ge \cdots \ge j_q >0$. For such a multi-index, we shall denote  $|J|:= j_1 + ...+ j_q$, $d(J):=q$, $\partial_s^{J}x := \partial_s^{j_1}x \dots \partial_s^{j_q}x$, $\ona^{J}x := \ona^{j_1}x\otimes \dots \otimes\ona^{j_q}x$.

Computation by induction, taking into account that $\ona_{\partial_s} \partial_s = k N$ and $\ona_{\partial_s} N = - k {\partial_s}$, gives that there are constants $c_{iJK}$ (different for each formula) such that
\begin{align}\lb{dsmp}
&\qquad \qquad \partial_s^m\psi = \ona^m\psi(\partial_s, ..., \partial_s) + \sum_{m,0}^{1,m-2} c_{iJK} k_\psi^i\ \partial_s^{J} k_\psi C(\ona^{K}\psi), 
\end{align}
\begin{align}
&\partial_s^m \(\ona^2\psi(N,N)\)  = \ona^{m+2} \psi(\partial_s, ..., \partial_s,N,N) +  \sum_{m,2}^{2,m-1} c_{iJK} k_\psi^i\ \partial_s^{J} k_\psi C(\ona^{K}\psi) \lb{dsmhp},
\end{align}
\begin{align}
\partial_s^m k &= \partial_s^m\(k_\psi+\<\ona\psi,N\>\)  = \partial_s^m k_\psi +\ona^{m+1} \psi(\partial_s, ..., \partial_s,N) +  \sum_{m,1}^{1,m-1} c_{iJK} k_\psi^i\ \partial_s^{J} k_\psi C(\ona^{K}\psi) \lb{dsmhpb},
\end{align}
\begin{align}\lb{dsmpb}
 \partial_s^m K = \ona^m K(\partial_s, ..., \partial_s) + \sum_{m,0}^{1,m-2} c_{iJK\ell} k_\psi^i\ \partial_s^{J} k_\psi C(\ona^{K}\psi \otimes\ona^{\ell}K) ,
\end{align}
where $C(\ona^{K}\psi)$ means $\ona^{K}\psi$ acting on $|K|$ copies of $\partial_s$ and/or $N$ and $\ds\sum_{m,r}^{s,t}$ means the sum over the indices satisfying $
i+|J|+d(J)+|K|+\ell=m+r$, $0\le d(J) \le [(m+r-s)/2]$,  $|K|+\ell\ge s$, $o(J)\le t$. Moreover, in this sums we consider that $\partial_s^Jk_\psi$ do not appear if $|J|=0$. 

For $d(K)$ one has $1\le d(K) \le m$, in \eqref{dsmp} and $1\le d(K) \le m+1$ in \eqref{dsmhp} and 	\eqref{dsmhpb}, and $0\le d(K) \le m-1$ in \eqref{dsmpb}.

This gives, for the derivatives of $k^2$, including the binomial coefficients in the constants, and renaming the constants $c_{iJK}$ in the last equality 
\begin{align}
\partial_s^m k^2 &= \sum_{\ell=0}^m \binom{m}{\ell} \partial_s^\ell k \partial_s^{m-\ell}k \nn \\
&= \sum_{\ell=0}^m \binom{m}{\ell} \(\partial_s^\ell k_\psi + C(\ona^{\ell+1}\psi) + \sum_{\ell,1}^{1,\ell-1}  c_{iJK} k_\psi^i\ \partial_s^{J} k_\psi C(\ona^{K}\psi) \) \nn \\
 & \qquad \qquad \qquad \( \partial_s^{m-\ell} k_\psi + C(\ona^{{m-\ell}+1}\psi) + \sum_{m-\ell,1}^{1,m-\ell-1}  c_{iJK} k_\psi^i\ \partial_s^{J} k_\psi C(\ona^{K}\psi) \) \nn \\
&= 2 \(C(\ona\psi) + k_\psi\) \partial_s^m k_\psi + \sum  c_{iJK} k_\psi^i\ \partial_s^{J} k_\psi C(\ona^{K}\psi) ,
\end{align}
where the last sum is along all the coefficients $i, J, K$ with $i\le m$, $|J| \le m-1$, $o(J)\le m-1$, and   some of the  $c_{iJK}$ are zero.

And, plugging all this into \eqref{ter1}, \eqref{ter2} and \eqref{ter3}, we obtain

\begin{align}
\partial_s^n (\Delta_\psi k_\psi)  &= \Delta_\psi (\partial_s^n k_\psi)  + n \(\ona^2\psi(\partial_s,\partial_s) +   k_\psi\ C(\ona\psi)+C(\ona \psi \otimes \ona \psi) \)\  \partial_s^{n} k_\psi  \nn \\
& \qquad \qquad \ +\sum_{\ell=2}^n \binom{n}{\ell} 
\(\ona^{1+\ell}\psi(\partial_s, ..., \partial_s) + \sum_{{1+\ell},0}^{1,\ell-1} c_{iJK} k_\psi^i\ \partial_s^{J} k_\psi C(\ona^{K}\psi) \) \partial_s^{n+1-\ell}k_\psi , \lb{ter1m}
\end{align}

\begin{align}
\partial_s^n \( k_\psi \ (k^2 - \oRic_\psi(N,N)) \)&= \(k_\psi^2+ 2 \<\ona\psi,N\> k_\psi + \<\ona\psi,N\>^2- K + \ona^2\psi(N,N)\) (\partial_s^n k_\psi) \nn \\
&+  2 \(C(\ona\psi) + k_\psi\) k_\psi \partial_s^n k_\psi + 2 \sum_{\ell=1}^{n-1} \binom{n}{\ell} \(C(\ona\psi) + k_\psi\) \partial_s^\ell k_\psi \partial_s^{n-\ell} k_\psi \nn \\
& + \sum_{\ell=1}^{n}\( \sum  c_{iJK} k_\psi^i\ \partial_s^{J} k_\psi C(\ona^{K}\psi)    \right.\nn
\\ &-  \(\ona^\ell K(\partial_s, ..., \partial_s) + \sum_{\ell,0}^{1,\ell-2} c_{iJKr} k_\psi^i\ \partial_s^{J} k_\psi C(\ona^{K}\psi \otimes\ona^{r}K)\) \nn\\
 & + \left. \(\ona^{\ell+2} \psi(\partial_s, ..., \partial_s,N,N) +  \sum_{\ell,2}^{2,\ell-1} c_{iJK} k_\psi^i\ \partial_s^{J} k_\psi C(\ona^{K}\psi)\) \) \partial_s^{n-\ell} k_\psi, \lb{ter2m}
\end{align}

\begin{align}
 \sum_{r=0}^{n-1} &\partial_s^r (k_\psi k \partial_s^{n-r}k_\psi) = k_\psi k \partial_s^nk_\psi + \sum_{r=1}^{n-1} \sum_{j=0}^r\binom{r}{j}\partial_s^j (k_\psi  \partial_s^{n-r}k_\psi)
 \partial_s^{r-j}k  \nn \\
 &= n (k_\psi^2+ k_\psi \<\ona\psi, N\>) \partial_s^nk_\psi + (k_\psi + \<\ona\psi, N\>) \sum_{r=1}^{n-1} \sum_{\ell=1}^r\binom{r}{\ell} \partial_s^{\ell}k_\psi \partial_s^{n-\ell}k_\psi \nn\\
 &+ \sum_{r=1}^{n-1} \sum_{j=0}^{r-1}\binom{r}{j}\partial_s^j (k_\psi  \partial_s^{n-r}k_\psi) \( \partial_s^{r-j} k_\psi +\ C(\ona^{{r-j}+1} \psi) +  \sum_{r-j,1}^{1,r-j-1} c_{iJK} k_\psi^i\ \partial_s^{J} k_\psi C(\ona^{K}\psi) \). \lb{ter3m}
\end{align}

We observe that, with the exception of the term $\Delta_\psi (\partial_s^n k_\psi)$, all these expressions \eqref{ter1m}, \eqref{ter2m} and \eqref{ter3m} have the form
\begin{align}
(a_n  + a_{n1} k_\psi +a_{n2} k_\psi^2) \partial_s^n k_\psi + \sum a_{iJ} k_\psi^i \partial_s^Jk_\psi,
\end{align}
where $i+|J|\ge 1$, $i\le n+1$, $o(J) \le n-1$, $|J| \le n$, and the coefficients \lq\lq$a_{\cdot \cdot}$'' are polynomials in the variables $\ona^m K$, $\ona^m \psi$, $m=1,...,n+1$  acting on $\partial_s$ and/or $N$, and  some of them can be zero. Then, using also \eqref{divtp2},  we can write (again renaming  the coefficients)
\begin{align}
\dfrac{d}{dt}\int_{\gamma(\cdot,t)}(\partial_{s}^n k_{\psi})^{2}ds_{\psi}&=
\int_{\gamma(\cdot,t)}\left[2 \partial_{s}^n k_{\psi}\(\Delta_\psi(\partial_s^nk_\psi)+ (a_n  + a_{n1} k_\psi +a_{n2} k_\psi^2) \partial_s^n k_\psi \right.\right. \nn \\
& \qquad \qquad +\left. \left.  \sum a_{iJ} k_\psi^i \partial_s^Jk_\psi\)  - k_\psi^2 (\partial_s^n k_\psi)^2   \right] ds_\psi \nn \\ 
&=
- 2 \int_{\gamma(\cdot,t)} \(\partial_{s}^{n+1} k_{\psi}\)^2 ds_\psi + 2 \int_{\gamma(\cdot,t)} (a_n  + a_{n1} k_\psi +a_{n2} k_\psi^2) (\partial_s^n k_\psi)^2 ds_\psi \nn \\ 
 & \qquad \qquad \qquad + \sum a_{iJ}  \int_{\gamma(\cdot,t)} \(k_\psi^i \partial_s^Jk_\psi\) (\partial_s^n k_\psi) ds_\psi ,
\label{E_igualdadprincipal}
\end{align}
with $i+|J|\ge 1$ in the last addend.
Now, we shall estimate each one of the addends in \eqref{E_igualdadprincipal}. First we observe that, if 
 $s_{0}\in\mathbb{S}^{1}$ is a critical point of $\partial_s^{j-1}k_{\psi}(\cdot)$. For any other point $s_1\in \mathbb{S}^{1}$, $j\ge 1$, we have:
\begin{align*}
\vert\partial^j_{s}(k_{\psi})(s_{1})\vert&=
\vert\partial^j_{s}(k_{\psi})(s_{1})-\partial^j_{s}(k_{\psi})(s_{0})\vert=
\Big\vert\int_{s_{0}}^{s_{1}}\partial_{s}^{j+1}k_{\psi}ds\Big\vert
\leq\int_{s_{0}}^{s_{1}}\vert\partial_{s}^{j+1}k_{\psi}\vert ds\leq
\dfrac{1}{E}\int_{s_{0}}^{s_{1}}\vert\partial_{s}^{j+1}k_{\psi}\vert e^{\psi} ds
\nn\\
&\leq \dfrac{1}{E}L_{\psi}(\gamma\vert_{[s_{0},s_{1}]})^{1/2}\Bigg(\int_{s_{0}}^{s_{1}}\vert\partial_{s}^{j+1}k_{\psi}\vert^{2} e^{\psi} ds\Bigg)^{1/2}
\leq \dfrac{1}{E}L_{\psi}(\gamma)^{1/2}\Bigg(\int_{\gamma}\vert\partial_{s}^{j+1}k_{\psi}\vert^{2} e^{\psi} ds\Bigg)^{1/2},
\end{align*}
where we have used the hypothesis  $0<E\leq e^{\psi}$ and Hölder's inequality. From this we conclude that  
\begin{align}
\max_{\mathbb{S}^{1}}(\partial_{s}^j(k_{\psi}))^{2}&\leq \dfrac{1}{E^{2}}L_{\psi}(\gamma)\int_{\gamma}(\partial_{s}^{j+1}k_{\psi})^{2}ds_{\psi} \qquad \text{ for } j\ge 1, \lb{desu}
\end{align}
and, using this bound, we obtain
\begin{align}\int_{\gamma}(\partial_{s}^j k_{\psi})^{2}ds_{\psi}&\leq
\dfrac{1}{E^{2}}L_{\psi}(\gamma)^{2}\int_{\gamma}(\partial_{s}^{j+1}k_{\psi})^{2}ds_{\psi}\qquad \text{ for } j\ge 1, \lb{desd}
\\
\int_{\gamma}(k_{\psi})^{2}(\partial_{s}^nk_{\psi})^{2}ds_{\psi}&\leq \dfrac{1}{E^{2}}L_{\psi}(\gamma)\int_{\gamma}k_{\psi}^{2}ds_{\psi}\int_{\gamma}(\partial_{s}^{n+1}k_{\psi})^{2}ds_{\psi}. \lb{dest}
\end{align}
We shall need also the following bound, alternative  to \eqref{desd}, that we shall obtain using  the bounds $0<E\le e^\psi\le D$ in \eqref{hipM1}, H\"older's inequality and integration by parts,
\begin{align}
\int_{\gamma}(\partial_{s}^nk_{\psi})^{2}ds_{\psi}&\leq
D\int_{\gamma}(\partial_{s}^n k_{\psi})^{2}ds\leq
D\Bigg(\int_{\gamma} (\partial_s^{n-1}k_{\psi})^{2}ds\Bigg)^{1/2}\Bigg(\int_{\gamma}(\partial^{n+1}_{s}k_{\psi})^{2}ds\Bigg)^{1/2}
\nn\\
&\leq
\dfrac{D}{E}\Bigg(\int_{\gamma}(\partial_s^{n-1}k_{\psi})^{2} ds_{\psi}\Bigg)^{1/2} \Bigg(\int_{\gamma}(\partial^{n+1}_{s}k_{\psi})^{2}ds_{\psi}\Bigg)^{1/2} 
\nn\\
&\leq
\dfrac{D}{E}\Bigg(\int_{\gamma}(\partial_s^{n-1}k_{\psi})^{2} ds_{\psi}\Bigg)^{1/2}\Bigg(1+\int_{\gamma}(\partial^{n+1}_{s}k_{\psi})^{2}ds_{\psi}\Bigg). \lb{L_desitres}
\end{align}

These inequalities will allow us to bound the second addend in \eqref{E_igualdadprincipal}. Now, to study the third addend, we use the inequalities \eqref{bmaxkp}, \eqref{desu}, \eqref{desd}  plus  Hölder's and Young's inequalities to obtain

\begin{align}
\sum& a_{iJ}
\int_{\gamma(\cdot,t)}k_{\psi}^{i}\partial_{s}^{J}k_{\psi}
\partial_{s}^{n}k_{\psi}ds_{\psi}
\leq \sum \dfrac{a_{iJ}}{2}
\Bigg\lbrace
\int_{\gamma(\cdot,t)}(\partial_{s}^{n}k_{\psi})^{2}ds_{\psi}
+\int_{\gamma(\cdot,t)}k_{\psi}^{2i}(\partial_{s}^{J}k_{\psi})^{2}ds_{\psi}
\Bigg\rbrace
\nn\\
&=\big(\sum \dfrac{a_{iJ}}{2}\big)\int_{\gamma(\cdot,t)}(\partial_{s}^{n}k_{\psi})^{2}ds_{\psi}
+\sum_{i\neq 0} \dfrac{a_{i0}}{2}\int_{\gamma(\cdot,t)}k_{\psi}^{2i}ds_{\psi}
+\sum_{d(J)\geq 1}\dfrac{a_{iJ}}{2} \int_{\gamma(\cdot,t)}k_{\psi}^{2i}(\partial_{s}^{J}k_{\psi})^{2}ds_{\psi}
\nn\\
&\leq
\big(\sum \dfrac{a_{iJ}}{2}\big)\int_{\gamma(\cdot,t)}(\partial_{s}^{n}k_{\psi})^{2}ds_{\psi}
+\sum_{i\neq 0} \dfrac{a_{i0}}{2}
\Big(c_{1}
+c_{2}\int_{\gamma(\cdot,t)}(\partial_{s}k_{\psi})^{2}ds_{\psi}
\Big)^{i-1}
\int_{\gamma(\cdot,t)}k_{\psi}^{2}ds_{\psi}
\nn\\
&\qquad +\sum_{d(J)\geq 1} \dfrac{a_{iJ}}{2}\dfrac{L_{\psi}(\gamma)^{q-1}}{E^{2(q-1)}}\Big(c_{1}
+c_{2}\int_{\gamma(\cdot,t)}(\partial_{s}k_{\psi})^{2}ds_{\psi}
\Big)^{i}  \nn \\ 
& \qquad \qquad \qquad \qquad \qquad \qquad \quad \times \int_{\gamma(\cdot,t)}(\partial_{s}^{j_{1}}k_{\psi})^2 ds_{\psi}
\int_{\gamma(\cdot,t)}(\partial_{s}^{j_{2}+1}k_{\psi})^2 ds_{\psi}
\cdots\int_{\gamma(\cdot,t)}(\partial_{s}^{j_{q}+1}k_{\psi})^2 ds_{\psi},\lb{67}
\end{align}
where we want to recall that  $d(J)\geq 1$ implies $n>1$ because $o(J)\le n-1$.

By substitution of the inequalities  \eqref{desd} to \eqref{67} in \eqref{E_igualdadprincipal} and taking into account that $k_\psi \le 1 + k_\psi^2$, we obtain

\begin{align}
\dfrac{d}{dt}&\int_{\gamma(\cdot,t)}(\partial_{s}^n k_{\psi})^{2}ds_{\psi}\leq
- 2 \int_{\gamma(\cdot,t)} \(\partial_{s}^{n+1} k_{\psi}\)^2 ds_\psi + 2(a_{n1} +a_{n2})\dfrac{L_{\psi}(\gamma)}{E^{2}}\int_{\gamma(\cdot,t)} k_\psi^2  ds_\psi
\int_{\gamma(\cdot,t)}(\partial_{s}^{n+1} k_\psi)^2 ds_{\psi}
\nn\\ 
&
\quad+\big(2a_{n}+2a_{n1}+\sum\dfrac{a_{iJ}}{2}\big)\dfrac{D}{E}\Bigg(\int_{\gamma(\cdot,t)}(\partial_{s}^{n-1}k_{\psi})^{2}ds_{\psi}\Bigg)^{1/2}
\Bigg(1+\int_{\gamma(\cdot,t)}(\partial_{s}^{n+1}k_{\psi})^{2}ds_{\psi}\Bigg)
\nn\\
&\quad+\sum_{i\neq 0} \dfrac{a_{i0}}{2}
\Big(c_{1}
+c_{2}\int_{\gamma(\cdot,t)}(\partial_{s}k_{\psi})^{2}ds_{\psi}
\Big)^{i-1}
\int_{\gamma(\cdot,t)}k_{\psi}^{2}ds_{\psi}
\nn\\
&\quad+\sum_{d(J)\geq 1} \dfrac{a_{iJ}}{2}\dfrac{L_{\psi}(\gamma)^{q-1}}{E^{2(q-1)}}\Big(c_{1}
+c_{2}\int_{\gamma(\cdot,t)}(\partial_{s}k_{\psi})^{2}ds_{\psi}
\Big)^{i} \nn \\
& \qquad \qquad \qquad \qquad \qquad \qquad\quad \times \int_{\gamma(\cdot,t)}(\partial_{s}^{j_{1}}k_{\psi})^2 ds_{\psi}
\int_{\gamma(\cdot,t)}(\partial_{s}^{j_{2}+1}k_{\psi})^2 ds_{\psi}
\cdots\int_{\gamma(\cdot,t)}(\partial_{s}^{j_{q}+1}k_{\psi})^2 ds_{\psi}. \lb{78}
\end{align}

When $n=2$, it appears a term with $j_2=1$ which gives rise to a $\int_{\gamma(\cdot,t)}(\partial_{s}^{n}k_{\psi})^2ds_\psi$ in the last line of the above expression. In this case we shall apply again \eqref{L_desitres} to this term. 
 From the induction hypothesis, given any $\eps_n>0$, there is $t^*>0$ such that, for every $t\ge t^*$ the sum of the coefficients of $\int_{\gamma(\cdot,t)} (\partial_s^{n+1}k_\psi)^2ds_\psi$ which are not in the first addend of \eqref{78} is lower than $1$,  and the sum of the terms which do not contain $\int_{\gamma(\cdot,t)} (\partial_s^{n+1}k_\psi)^2ds_\psi$  nor $\(\int_{\gamma(\cdot,t)}(\partial_s^{n-1}k_{\psi})^{2}ds_{\psi}\)^{1/2}$, is also lower than $\eps_n$. Using that, we can write \eqref{78} as 
\begin{align}\lb{E_clavekpsiceroj}
\dfrac{d}{dt}\int_{\gamma(\cdot,t)}(\partial_{s}^n k_{\psi})^{2}ds_{\psi}&\leq
\eps_n + C_n \(\int_{\gamma(\cdot,t)}(\partial_s^{n-1}k_{\psi})^{2}ds_{\psi}\)^{1/2} ,
\end{align}
with $C_n>0$. On the other hand, if we repeat the above process to obtain  \eqref{78}, but without using \eqref{L_desitres} and without forgetting the negative term containing the derivative of highest order, we obtain
\begin{align}\lb{E_clavederivadakpsicero}
\dfrac{d}{dt}\int_{\gamma(\cdot,t)}(\partial_{s}^{n-1} k_{\psi})^{2}ds_{\psi}&\leq - \int_{\gamma(\cdot,t)}(\partial_s^{n}k_{\psi})^{2} ds_\psi +
\eps_{n-1} + D_n \int_{\gamma(\cdot,t)}(\partial_s^{n-1}k_{\psi})^{2}ds_{\psi}.
\end{align}

If we denote 
$\ds g_{j}(t):=\int_{\gamma(\cdot,t)}(\partial_s^jk_{\psi})^{2}ds_{\psi}$, $j=0,1, ..., n, ... $,  
the inequalities (\ref{E_clavekpsiceroj}) and (\ref{E_clavederivadakpsicero}) become\begin{align}
\label{E_clavekpsiceronot}
g_{n-1}'(t)&\leq -g_n(t)+ D_n g_{n-1}(t) + \eps_{n-1},
\\
\label{E_clavederivadakpsiceronot}
g_n'(t)&\leq \eps_n + C_n g_{n-1}(t)^{1/2}.t
\end{align}
Using again that $\lim_{t\rightarrow\infty}g_{n-1}(t)=0$, for every $\eps>0$ there is a $t_{0}=t_{0}(\eps)$ such that $0\leq g_{n-1}(t)\leq\eps^{2}$, for every $t\geq t_{0}$. By (\ref{E_clavederivadakpsiceronot}), $g_n'(t)\leq C_n \eps + \eps_n$ for every $t\geq t_{0}$. By integration we obtain
\begin{equation}\label{E_convergenciapsigeodesicasvariag}
g_n(t_{2})-g_n(t_{1})\leq (C_n \eps + \eps_n) (t_{2}-t_{1}) ,\:\text{ for every }\: t_{2}\geq t_{1}\geq t_{0}
\end{equation}
and, from (\ref{E_clavekpsiceronot}):
\begin{equation}\label{E_convergenciapsigeodesicasvariaf}
g_{n-1}'(t)\leq-g_n(t)+ D_n\eps^{2} + \eps_{n-1},\:\text{ for every } t \geq t_{0}.
\end{equation}

Now we claim that $\lim_{t\rightarrow\infty} g_n(t)=0$. We shall prove it by contradiction. Let us suppose that there are $C>0$ and $\lbrace t_{k}\rbrace_{k=1}^{\infty}$, $t_{k+1}\geq t_{k}$, $t_{1}\geq t_{0}$ such that $\lim_{k\rightarrow\infty}t_{k}=\infty$ and $g_n(t_{k})\geq C$, for every $ k\in\mathbb{N}$.
From (\ref{E_convergenciapsigeodesicasvariag}) we have
\begin{equation*}
g_n(t)\geq\dfrac{C}{2},\quad \text{ for every } t\in[max\lbrace t_{0},t_{j}-\dfrac{C}{2(C_n \eps+\eps_n)}\rbrace,t_{j}]\neq\emptyset,\:\text{ for every } j\in\mathbb{N}
\end{equation*}
Then, from (\ref{E_convergenciapsigeodesicasvariaf}),
\begin{align}
g_{n-1}'(t)&\leq-\dfrac{C}{2} + D_n\eps^{2} + \eps_{n-1},\lb{thisine} \\ 
\text{ for every }&t\in[max\lbrace t_{0},t_{j}-\dfrac{C}{(C_n \eps+\eps_n)}\rbrace,t_{j}]\neq\emptyset,\:\text{ and } j\in\mathbb{N}. \nn
\end{align}
Let us choose $\eps$, $\eps_{n-1}$ and $j(\eps)$ such that  $A(\eps):=\dfrac{C}{2}- (D_n\eps^{2} + \eps_{n-1})>0$ and $s_{j}:=t_{j}-\dfrac{C}{2(C_n \eps+\eps_n)}> t_{0}$ for every $j\ge j(\eps)$. Then $g_{n-1}'(t)\leq -A(\eps)$ and $g_{n-1}(t) \le \eps^2$ for every $t\in[s_{j}(\eps),t_{j}]$, $j\ge j(\eps)$. By integration of the inequality \eqref{thisine}
\begin{equation}\lb{ttj0}
0\leq g_{n-1}(t)\leq g_{n-1}(s_{j}(\eps))-A(\eps)(t-s_{j}(\eps)),\:\text{ for every }\:t\in[s_{j}(\eps),t_{j}].
\end{equation}
 For $\eps$ small enough, the last term of \eqref{ttj0} vanishes when 
$\ds t=s_{j}(\eps)+\dfrac{g_{n-1}(s_{j}(\eps))}{A(\eps)} \le s_{j}(\eps)+\dfrac{\eps^2}{A(\eps)} < t_j$. For these values of $\eps$ and $j(\eps)$, it follows again from \eqref{ttj0} that, for $j\ge j(\eps)$
\begin{align*}
0\leq g_{n-1}(t)\leq g_{n-1}(s_{j}(\eps))-A(\eps)(t-s_{j}(\eps))<0,\:\text{ for every }\:t\in ]s_j(\eps)+ \fracc{\eps^2}{A(\eps)},t_j[
\end{align*}
which is a contradiction. Therefore $\lim_{t\rightarrow\infty} g_n(t)=0$, which finishes the proof of Step (3) by induction.
\qed

{\bf Step 4.} \ {\it Let  $\gamma:\mathbb{S}^{1}\times[0,\infty)\rightarrow \oM^{2}$ be a solution of the \DMCF satisfying (\ref{ProSol2}) and (\ref{ProSol4}) moving in a surface with density $(\oM^{2},g,e^{\psi})$  satisfaying \eqref{hipM1}. 
If this solution exists for every $t\in[0,\infty[$, then, for every $m=0,1,2,...$,
$\partial_s^mk_{\psi}$ converges uniformly to zero when  $t\rightarrow\infty$.
}
 
 In fact, this is a consequence of  \eqref{P_kpsicef} and \eqref{acotg}, the decreasing property of $L_\psi$, the condition \eqref{ProSol2}  and the following inequality whose proof follows the ideas used to prove \eqref{bmaxkp}. For every $t$, let $s_0$ be the point where $\ds |\partial_s^m k_\psi(s_0,t)|=\min_{s\in \mathbb{S}^{1}}|\partial_s^mk_\psi(s,t)|$, one has
\begin{align*}
\vert \partial_s^m k_{\psi}(s_{1})\vert &= 
\Big\vert \partial_s^m k_{\psi}(s_{0}) +\int_{s_{0}}^{s_{1}}\partial_{s}^{m+1} k_{\psi}ds_{t}\Big\vert \le  \frac1{L_\psi}\int_{\gamma(\cdot,t)}|\partial_s^mk_\psi| ds_{\psi} + \frac1{E}\int_{\gamma(\cdot,t)} |\partial_s^{m+1}k_\psi| ds_\psi \\
& \le \(\frac1{L_\psi}\int_{\gamma(\cdot,t)} (\partial_s^m k_\psi)^2 \ ds_\psi\)^{1/2} + \frac1E  \(L_\psi \int_{\gamma(\cdot,t)} (\partial_s^{m+1}k_\psi)^2 ds_\psi\)^{1/2}
\end{align*}
 for every point $\gamma(s_1,t)$ in $\gamma(\cdot,t)$. By Step 3, this goes to zero when $t\to\infty$, and step 4 is proved.  
\qed

{\bf Step 5: End of proof of Theorem \ref{ctapg}.} From the hypothesis that the evolution of $\gamma$ is contained in a compact subset of $\oM$ it follows that $\gamma(\cdot,t)$ is uniformly bounded. To show that $|\partial_\theta\gamma(\theta,t)|$ is also bounded, we observe that  $ds_\psi=|\partial_\theta\gamma(\theta,t)| e^\psi d\theta$ and the variation formula \eqref{evol.dV} gives $\partial_t (|\partial_\theta\gamma(\theta,t)| e^\psi) = - k_\psi^2 |\partial_\theta\gamma(\theta,t)| e^\psi$, then $\partial_t \ln(|\partial_\theta\gamma(\theta,t)| e^\psi) =  -  k_\psi^2$. Therefore $|\partial_\theta\gamma(\theta,t)| e^\psi$ is  decreasing with $t$  and, since $e^\psi$ is bounded by hypothesis, 
\begin{align}
|\partial_\theta\gamma(\theta,t)| \text{ is also bounded by a constant independent of $t$}. \lb{tam}
\end{align}

We now reparametrize $\gamma$ by $u=\fracc{s_t}{2 \pi L_t}$. That is, we consider $\wga(u,t)= \gamma(\p_t(u),t)$, where, for every $t\ge 0$, $\p_t$ is the inverse of the function $\theta\mapsto \fracc{\int_0^\t|\partial_\theta\gamma(\theta,t)| d\t}{2 \pi \int_0^{2 \pi}|\partial_\theta\gamma(\theta,t)| d\t} = \fracc{s_t(\theta)}{2 \pi L_t}$. Since $\wga(\cdot,t)$ and $\gamma(\cdot, t)$ are geometrically the same curve, all the estimates we have for $k$, $k_\psi$, $\partial^m_s k$ and $\partial_s^m k_\psi$ are the same for both curves. Moreover, 
\begin{align}
\parcial{s_t}{u} &= 2 \pi L_t \quad \text{ , } \quad \parcial{^m s_t}{u}=0 \text{ for } m\ge 2, \\
\partial_{s_t}^2\gamma &= k N, \\
\partial_{s_t}^3\gamma &= \partial_sk\ N - k\ T, \\
\partial_{s_t}^4\gamma &= \partial_s^2 k\ N -  \partial_sk\ k\ T - \partial_sk\ T - k^2 \ N,\\
\partial_{s_t}^5\gamma &= \partial_s^3 k\ N -2\ \partial_s^2 k\ k \ T - (\partial_s k)^2 T - \partial_s k\ k^2 \ N - \partial_s^2 k \ T - 3 k \ \partial_sk\ N + k^3 T,\\
\dots & \dots \dots \nn \\
\partial_{s_t}^m\gamma &= \xi_m(k, \partial_s k, ..., \partial_s^{m-3}k) T + \zeta_m (k, \partial_s k, ..., \partial_s^{m-2}k) N, \\
\dots & \dots \dots \ .\nn
\end{align}
where $\xi_m$ and $\zeta_m$ are polynomials in $k, \partial_s k, ..., \partial_s^{m-2}k$ of degree lower than $m-2$, where the degree of each monomial is obtained counting the degree of $\partial_s^jk$ as $j+1$.
Using these formulae for the computation of $\partial_u\wga$, we obtain  
\begin{align}
|\partial_u\wga| &= |\parcial{s_t}{u}\partial_{s_t}\gamma | = |2 \pi L_t \partial_{s_t}\gamma | = 2 \pi L_t , \lb{cure}\\
|\partial^2_u\wga| &= | (2 \pi L_t)^2 \partial_{s_t}^2\gamma | = (2 \pi L_t)^2 |k| , \nn \\
\dots & \dots \dots \nn \\
|\partial_u^m\wga| &= (2 \pi L_t)^m \sqrt{\xi_m^2+\zeta_m^2}.
\end{align}

From Step 4 and \eqref{dsmhpb} follow that $\sqrt{\xi_m^2+\zeta_m^2}$ is bounded by a bound independent of $t$. From \eqref{hipM1} we have $L_t= \int_{\gamma(\cdot,t)} e^{-\psi} ds_\psi \le \fracc1{E}L_\psi(\gamma(\cdot,0))$. Then, for every $m$, by Ascoli-Arzela, there is a sequence $\wga(\cdot,t_n)$ which converges to a $C^m$ curve that, by Step 4, has $k_\psi=0$. Moreover, by \eqref{cure} and \eqref{ProSol2}, the limit curve is regular. This finishes the proof of Theorem \ref{ctapg}.


\section{Proof of theorems \ref{plnr} and \ref{pgen}}\lb{sec4}
An important tool in the proof of theorem 1 is the variation of the area enclosed by a curve which evolves under \DMCF. In this case,  formula \eqref{igu_vol} becomes
\begin{align}\lb{var_V_r}
\deri{A_\og}{t} = - \int_\gamma k_\psi ds = -\int_\gamma k\ ds +\int_\gamma \<\ona\psi, N\> ds = -2\pi - \int_{\Omega_t} \oDelta \psi da_\og,
\end{align}
where we have used the divergence theorem and the differentiability of $\psi$ on $\Omega_t$ in the last equality (with the sign "-" because $N$ points inward). If $\psi$ is radial, $\oDelta \psi= \tr(\ona^2\psi) = \tr(\ona (\psi'\ona r)) = tr (\psi'' \ona r\otimes\ona r+ \psi'\ona^2r) ) = \psi'' +\psi' \fracc{1}{r}$, then substitution in \eqref{var_V_r} gives 
 \begin{align}\lb{var_V_rf}
\deri{A_\og}{t} =-2\pi - \int_{\Omega_t} \(\psi''+ \frac1r \psi'\) da_\og.
\end{align}
Looking for densities with a nice variation of $A_\og(\Omega_t)$, we can consider the solutions of $\psi''+ \fracc1r \psi'=\lambda$ (where $\lambda$ is a constant), which give $\deri{A_\og}{t} =  -2\pi - \lambda A_\og(\Omega_t)$  and are 
\begin{equation}\lb{dAct}\psi (r) = \lambda \fracc{r^2}{4} + b + a \ln(r).\end{equation}
But this solution has a singularity at $r=0$ if $a\ne 0$, so that formula \eqref{var_V_r} cannot be applied, because it has been obtained using the differentiability of $\psi$ on all $\Omega_t$. For the case where we have a singularity at the origin, we need to distinguish the following situations:

S1) The origin is in the interior of $\Omega_t$, \eqref{var_V_r} has to be written in the following way:
\begin{align}\lb{var_V_s}
\deri{A_\og}{t} =  -\int_\gamma k\ ds +\int_\gamma \<\ona\psi, N\> ds = -2\pi - \lim_{\rho\to 0}\int_{\Omega_t-B^2(\rho)} \oDelta \psi da_\og - \lim_{\rho\to 0} \int_{\mathbb{S}^{1}(\rho)} \<\ona\psi, \xi\> ds,
\end{align}
where $\xi=\ona r$ is the unit vector field normal to the circle  $\mathbb{S}^{1}(\rho)$ of radius $\rho$ pointing into the interior of $\Omega_t-B^2(\rho)$. If $\psi$ is radial, this gives:
\begin{align}
\deri{A_\og}{t} &=    -2\pi - \lim_{\rho\to 0}\int_{\Omega_t-B^2(\rho)} \(\psi''+ \frac1r \psi'\)da_\og - \lim_{\rho\to 0} \int_{\mathbb{S}^{1}(\rho)} \psi' ds  \nn\\
&=-2\pi - \lim_{\rho\to 0}\int_{\Omega_t-B^2(\rho)} \(\psi''+ \frac1r \psi'\) da_\og - \lim_{\rho\to 0} 2 \pi \rho  \psi'(\rho), \lb{var_V_rs}
\end{align}
and, when $\psi$ has the form \eqref{dAct},
\begin{align}
\deri{A_\og}{t} & = - 2\pi - \lambda A_\og(\Omega) -  2 \pi \ a . \lb{var_V_rsct}
\end{align}

S2) The origin is outside $\Omega_t$, then the right formula is \eqref{var_V_r} and, when $\psi$ has the form 
\eqref{dAct},
\begin{align}
\deri{A_\og}{t} & = - 2\pi - \lambda A_\og(\Omega) . \lb{var_V_rsc2}
\end{align}
Now, let us consider the cases of Theorem \ref{plnr}:

Case (1) $\lambda=0$, $a=-1$. 

In this case $\psi$ has a singularity at the origin. But given any curve contained in $\re^2-\{0\}$, by Proposition \ref{critmin} there is a $\psi$-minimal curve (bounding a disk $D$) between the origin and the curve which acts as a barrier, and the boundary of a disc containing the curve acts as another barrier. Then the evolution of the curve will be in the domain bounded by these two $\psi$-minimal curves, and  the hypotheses of theorems \ref{ThAO} and \ref{ctapg} are satisfied.

If the domain bounded by $\gamma_0$ contains the origin, the formula \eqref{var_V_rsct} says that the area of the domain $\Omega_t$ is constant, then the flow is defined for all time. By Theorem \ref{ctapg}, it must subconverge to a $\psi$-minimal curve enclosing the same area which, by Proposition \ref{critmin},  is the circle of radius $\sqrt{A/\pi}$.

 If the domain bounded by $\gamma_0$ does not contain the origin, then formula \eqref{var_V_rsc2} says that the area of the domain $\Omega_t$ decreases at constant velocity $2 \pi$, so, by Theorem \ref{ThAO}, the flow shrinks to a point in time $T=A/(2\pi)$.
 
 Case (2) $\lambda>0$, $a<-1$. 
 
Again $\psi$ has a singularity at the origin. By Proposition \ref{cilindros}, the circle of radius $r=\sqrt{-2 (a+1)/\lambda}$ is the unique circle $\psi$-minimal curve, and it is an attractor. Therefore, given any curve contained in $\re^2-\{0\}$, there is a circle of radius lower than $r_{min}$ and another of radius bigger than the maximum between $r_{max}$ and $\sqrt{-2 (a+1)/\lambda}$, which act as  barriers. Then  the hypotheses of theorems \ref{ThAO} and \ref{ctapg} are satisfied.

If the domain bounded by $\gamma_0$ contains the origin, the barriers indicated in the previous paragraph push $\gamma$ to the unique circle $\psi$-minimal curve, then, it cannot go to a point and, by Theorem \ref{ctapg}, it must subconverge to a closed $\psi$-minimal curve, that must be the mentioned circle.

 If the domain bounded by $\gamma_0$ does not contain the origin, the  formula \eqref{var_V_rsc2} gives a negative upper bound for the speed of the area. Then the flow exists only for finite time and Theorem \ref{ThAO} implies that the limit is a point. By \eqref{var_V_rsc2} $T_{\max}=\dfrac{1}{\lambda}ln\Big(1+\dfrac{\lambda A_{\overline{g}}(0)}{2\pi}\Big)$.  This finishes the proof of Theorem \ref{plnr}.

Let us start with the proof of Theorem \ref{pgen}.

From the general hypothesis of the Theorem and Proposition \ref{cilindros}, when $\psi$ is smooth, given any simple closed curve $\gamma_0$, there is always a circle bounding a disk that contains $\gamma_0$ such that, or it is a $\psi$-minimal curve, or it shrinks under \DMCF. Then $\gamma$ moves in a bounded domain, and we can apply theorems \ref{ThAO} and \ref{ctapg}. When $\psi$ is not smooth and $\psi'(t)<-1/t$ for $t\in]0,r_1[$, we have two circles bounding an annulus that contains  $\gamma_0$, which are $\psi$-minimal curves, or move up to a $\psi$-minimal curve and the situation  is similar to the smooth case.

Now we study in detail the cases where $\psi$ is smooth. 

If $r_{max} < r_1$, the maximal time of existence of the flow of $\gamma$ is bounded by the corresponding time for a circle of radius $r_{max}+\delta<r_1$, which moves with speed given by \eqref{drdtS}. Since $\psi'$ is continuous on $[0,r_{\max}+\delta]$, and $\psi'(r)+1/r>0$, there is an $\eps>0$ satisfying $\psi'+1/r>\eps$, so that the radius of the circle evolves according to $\deri{r}{t}=-(\frac1r+\psi) <-\eps$. Then the circle collapses to a point in finite time and, by the avoidance principle, the same happens with $\gamma$.

If $r_{max}\le r_1$ and equal to $r_1$ at some points of $\gamma$, we can use \eqref{evol.re} restricted to our situation with $n=1$ and $\psi$ radial to write
\bec\lb{drtcu}
\parcial{r}{t} = \Delta r + \frac{|\nabla r|^2}{r} - \< \ona r, N\>^2 \psi' - \frac1r. \eec
 To apply the maximum principle, we take into account that, in a point with maximal $r$, $\Delta r\le 0$, $\nabla r=0$ and $N=- \ona r$, which gives 
 $$\parcial{r}{t} \le  - \( \psi' + \frac1r\)\le 0 \text{ as far as } r\le r_1. $$
 Since at time $t=0$, $r\le r_1$, by the maximum principle $r(t) \le r_1$ and there are points with $r<r_1$. By the strong maximum principle, $r(t)<r_1$ for $t>0$, so we can apply the argument of the previous case staring with $\gamma(\cdot,t)$.

If $r_{2k-1}< r <r_{2k+1}$ and $0\notin \Omega_0$,  the curve is inside by the ring bounded by two circles of radius $r_{2k-1}+\delta$ and $r_{2k+1}-\delta$ that move to the circle of radius $r_{2k}$. There is a time such that the area of this ring is as small as we want. By \eqref{var_V_r}, the area enclosed by the curve decreases with finite speed, then the flow of $\gamma$ exists only for finite time and, by Theorem \ref{ThAO}, it collapses to a point.

If $r_{2k-1}\le r \le r_{2k+1}$ where at least one equality is satisfied at some point, we can use \eqref{drtcu} and the strong maximum priciple as above to conclude that, for $t>0$, $r_{2k-1} < r < r_{2k+1}$, and we apply  the argument of the previous paragraph.

If $r_1 < r_{min}$ and $0\in \Omega$, the circle $r=r_1$ and a circle of radius $r>r_{max}$ act as barriers. So $\gamma$ cannot collapse to a point and, by Theorem  \ref{ctapg}, it subconverges to a $\psi$-minimal curve.

If $r_{2k-1}< r <r_{2k+1}$ and $0\in \Omega_0$, the unique closed simple $\psi$-minimal curve contained in this ring is the circle of radius $r_{2k}$, which gives case iii.1.

If in the last two cases we have \lq\lq$\le$'' instead of \lq\lq$<$'', we argue as before using \eqref{drtcu} and the strong maximum principle to reduce the situation to the previous one.

If $\lim_{t\to 0}\psi'(t) = -\infty$ and $\psi'>-\frac1r$ for $r<r_1$, we cannot assure that a circle of radius $<r_1$ goes to $0$ in finite time. Also we cannot assure that $\gamma$ will move in a region with $\psi'$ bounded in order to apply theorems \ref{ThAO} and \ref{ctapg}. For this reason, we do not consider the situation $r_{max}\le r_1$ when $\psi$ has a singularity at the origin and $\psi'>-\frac1r$ for $r\in ]0,r_1]$. For the other cases the discussion is exactly the same that for the smooth case.

If $\lim_{t\to 0}\psi'(t) = -\infty$ and $\psi'<-\frac1r$ for $r<r_1$, the situation is the same that in the smooth case for curves with $r_{min}\ge r_1$, and this gives case b.ii.3). Cases b.ii.1) and  b.ii.2), respectively, by the same arguments that cases ii) and iii.1) in the smooth setting.
\qed

\section{Appendix: another view to (anti-)Gaussian density}

In $\mathbb{R}^{n+1}$ the (anti)-Gaussian density corresponds to $\psi(x)=\eps \frac12 n\mu^2 |x|^2$, $\eps=1$ anti-Gaussian and $\eps=-1$ Gaussian. These flows have two very particular properties which allowed to A. Borisenko and the first author to describe  the evolution of a compact convex hypersurface of $\re^{n+1}$ under these flows. The first particular property is

\begin{lema}\lb{1}{\rm(\cite{sm,bomi6})} $F(\cdot,t)$ is a solution of \eqref{gmcf} (with $\psi$ an (anti)-Gaussian density) iff $\hF(\cdot,\htt) = e^{\eps n\mu^2 t(\htt)} F(\cdot, t(\htt))$ is a solution of \begin{equation}\label{tmcf}
\<\parcial{\widehat{F}}{\widehat{t}},\widehat{N}\> = \widehat{H},
\end{equation}
 where $t(\htt) = \fracc1{\eps 2n  \mu^2} \ln( 1 + \eps 2n  \mu^2 \ \htt)$. (When $\eps=-1$, $\htt<\fracc1{2n\mu^2}$).
\end{lema}

\noindent The statement of Lemma \ref{1} given here is that of \cite{bomi6}, but it is based on a more general result on the equivalence of flows given by K. Smoczyk in \cite{sm}. The other special property is 

\begin{lema}\lb{2}{\rm(\cite{bomi6})}
 Let $p_0$ be a point  in $\re^{n+1}$.  Then the motion $F_t(M)$ of $F_0(M)$ is the composition of the motion $\widetilde F_t(M)$ of $F_0(M)-p_0$ by the flow  \eqref{gmcf} with the translation of vector $p(t) = e^{ - \eps n \mu^2 t} p_0$, that is $F_t(x)= e^{-\eps n \mu^2t} p_0 + \widetilde F_t(x)$. 
\end{lema}

 Here, we shall use again these properties to study the evolution of a closed embedded curve in the plane $\re^2$ under \DMCF, with $\psi$ an (anti)-Gaussian density. We shall prove:

\begin{teor}\lb{evpc} In the plane $\re^2$ with (anti)-Gaussian density $\psi=\eps \frac12 \mu^2 |x|^2$,
let $\gamma(\cdot,t)$ be a solution of \eqref{gmcf} such that $\gamma(\cdot,0)$ is an embedded curve. Let us denote by $A$ the area of the region bounded by $\gamma(\cdot,0)$. 
\begin{enumerate}
\item If $\eps=1$,  the maximal solution $\gamma(\cdot,t)$ is defined on $\mathbb{S}^{1}\times[0,T[$, \linebreak $T=\fracc1{ 2  \mu^2} \ln( 1 + \mu^2 \ \fracc{A}{\pi})$ and $\lim_{t\to T}\gamma(\mathbb{S}^{1},t)$ is a round point.
\item If $\eps=-1$:
\begin{enumerate}
\item If $A<\pi/\mu^2$,  the maximal solution $\gamma(\cdot,t)$ is defined on $\mathbb{S}^{1}\times[0,T[$, \linebreak $T=-\fracc1{ 2  \mu^2} \ln( 1 - \mu^2 \ \fracc{A}{\pi})$ and $\lim_{t\to T}\gamma(\mathbb{S}^{1},t)$ is a round point.
\item  If $A=\pi/\mu^2$,   the maximal solution $\gamma(\cdot,t)$ is defined on $\mathbb{S}^{1}\times[0,\infty[$, and $\lim_{t\to \infty}\gamma(\mathbb{S}^{1},t)$ is a circle of radius $1/\mu$ centered at the origin $(0,0)\in\re^2$ or centered at the infinite.
\item If $A>\pi/\mu^2$,   the maximal solution $\gamma(\cdot,t)$ is defined on $\mathbb{S}^{1}\times[0,\infty[$, and $\lim_{t\to \infty}\gamma(\mathbb{S}^{1},t)$ can be a straight line through $(0,0)\in\re^2$ or a curve in the infinite bounding a region in the infinite, or a curve in the infinite bounding the whole space.
\end{enumerate}
\end{enumerate}
\end{teor}
\begin{demo}

Let $A(t)$ (resp. $\hA(\htt)$) the area of the domain bounded by the curve $\gamma_t(\mathbb{S}^{1})$ (resp. $\widehat{\gamma}_\htt(\mathbb{S}^{1})$).  From the works \cite{gh,gr87,gr89}, it is known that the flow $\widehat{\gamma}_\htt$ of $\widehat{\gamma}_0$ under \eqref{tmcf} evolves to a round point when $\htt \to \fracc{\hA(0)}{2\pi}$. From the definition of $\widehat{\gamma}$ and $t(\htt)$ one has $\widehat{\gamma}(\cdot,0) = \gamma(\cdot,0)$, then $\widehat{A}(0)= A(0)$. Using Lemma \ref{1}, the evolution of $\gamma$ follows from the evolution of $\widehat{\gamma}$ according to the following cases:

When $\eps=1$, this gives that $\gamma_{t}$ is well defined for $t\in[0,\fracc1{2  \mu^2} \ln( 1 +  \mu^2 \ \fracc{A(0)}{\pi})[$, and the limit of $\gamma_t$ when $t\to \fracc1{ 2  \mu^2} \ln( 1 +   \mu^2 \ \fracc{A(0)}{\pi})$ is a round point.

When $\eps = -1$, the bound $\htt<\fracc1{2\mu^2}$ forces us to distinguish three cases:

1) $A(0) < \fracc{\pi}{\mu^2}$. The evolution is like in case $\eps=1$.

2) $A(0) = \fracc{\pi}{\mu^2}$. The solution $\gamma(\cdot,t)$ is well defined for $t\in[0,\infty[$. When the $\lim_{\htt\to{1/2\mu^2}} \widehat{\gamma}(\mathbb{S}^{1},\htt)=(0,0)$, we are just in the case $\gamma$ is the normalized motion associated to the mean curvature motion $\widehat{\gamma}$, as described in \cite{hu84} (see the last section of \cite{bomi6}), then  $\lim_{t\to\infty} \gamma(\mathbb{S}^{1},t)$ is a circle of radius $1/\mu$ centred at $(0,0)$. When $\lim_{\htt\to{1/2\mu^2}} \widehat{\gamma}(\mathbb{S}^{1},\htt)\ne (0,0)$, the above result combined with  Lemma \ref{2}  gives that $\lim_{t\to\infty}\gamma(\mathbb{S}^{1},t)$ is a circle of radius $1/\mu$ and center a point in the infinite.

3) $A(0)>\fracc{\pi}{\mu^2}$. Then $\lim_{\htt\to{1/2\mu^2}} \widehat{\gamma}(\mathbb{S}^{1},\htt)=:C_{1/2\mu^2}$ is a closed embedded curve (because the flow of $\gamma$ stops before the flow of $\widehat{\gamma}$ finishes). Let $\hOm_{1/2\mu^2}$ be the domain bounded by $\hC_{1/2\mu^2}$, and $\Omega_t$ the domain bounded by $\gamma(\mathbb{S}^{1},t)$. There are three possibilities:

\begin{description}
\item{3.1} $(0,0)\in C_{1/2\mu^2}$, then $\lim_{t\to\infty}\gamma(\mathbb{S}^{1},t)$ is a line and $\lim_{t\to\infty}\Omega_t$ is a half space.
\item{3.2} $(0,0)\in \Omega_{1/2\mu^2}- C_{1/2\mu^2}$, then $\lim_{t\to\infty}\gamma(\mathbb{S}^{1},t)$ is a curve in the infinite and $\lim_{t\to\infty}\Omega_t$ is the whole space $\re^2$.
\item{3.3} $(0,0)\notin \Omega_{1/2\mu^2}$, then $\lim_{t\to\infty}\gamma(\mathbb{S}^{1},t)$ is a curve in the infinite and $\lim_{t\to\infty}\Omega_t$ is also in the infinite.
\end{description}
\end{demo}


\bibliographystyle{alpha}

\noindent
Department of Geometry and Topology \\
University of Valencia \\
46100-Burjassot (Valencia), Spain\\
 {miquel@uv.es} and {Francisco.Vinado@uv.es}

\end{document}